%% file: MSML2022_arXiv.tex
\documentclass[onefignum,onetabnum]{siamart190516}

\input{ex_shared}

\ifpdf
\hypersetup{
  pdftitle={Online Weak-form Sparse Identification of Partial Differential Equations},
  pdfauthor={D. Messenger, E. Dall'Anese, D. Bortz}
}
\fi

\begin{document}

\maketitle

\begin{abstract}
This paper presents an online algorithm for  identification of partial differential equations (PDEs) based on the weak-form sparse identification of nonlinear dynamics algorithm (WSINDy). The algorithm is online in a sense that if performs the identification task by processing solution snapshots that arrive sequentially. The core of the method combines a weak-form discretization of candidate PDEs with an online proximal gradient descent approach to the sparse regression problem. In particular, we do not regularize the $\ell_0$-pseudo-norm, instead finding that directly applying its proximal operator (which corresponds to a hard thresholding) leads to efficient online system identification from noisy data. We demonstrate the success of the method on the Kuramoto-Sivashinsky equation, the nonlinear wave equation with time-varying wavespeed, and the linear wave equation, in one, two, and three spatial dimensions, respectively. In particular, our examples show that the method is capable of identifying and tracking systems with coefficients that vary abruptly in time, and offers a streaming alternative to problems in higher dimensions.  
\end{abstract}

\section{Context and Motivations}
\textcolor{white}{....} System identification (SID) and parameter estimation of dynamical systems are ubiquitous tasks in scientific research and engineering, and are required steps in many control frameworks. A typical strategy is to solve a regression problem based on sample trajectories from the underlying system, with few samples available in practice. Identification of dynamical systems is a classical field of research~\cite{ljung1999system}; recently, several  works provided new theoretical insights on the efficacy of classical first-order optimization methods in solving SID problems based on single trajectories  (see, e.g., \cite{fattahi2019learning,foster2020learning,sattar2020non,simchowitz2018learning} and references therein). Existing results in this context are heavily focused on discrete-time, finite-dimensional systems of known functional form, yet the focus on single-trajectory data paves the way for identification of more complex dynamical systems in the {\it online} setting, which is the subject of the current article.

By suitably discretizing candidate dynamical systems using data and employing sparse regression, SID and parameter estimation can be accomplished simultaneously. A notable development in this pursuit is the sparse identification of nonlinear dynamics (SINDy) algorithm (\cite{brunton2016discovering}), a general framework for discovering dynamical systems using sparse regression. Since its inception in the context of autonomous ordinary differential equations (ODEs), SINDy has been extended to autonomous partial differential equations (PDEs) (\cite{rudy2017data,schaeffer2017learning}), stochastic differential equations (SDEs) (\cite{boninsegna2018sparse}), non-autonomous systems (\cite{rudy2019data}), and coarse-grained equations (\cite{bakarji2021data}), to name a few.

A significant challenge in using SINDy to solve real-world problems is the computation of derivatives from noisy data. Within the last few years, the consensus has emerged that weak-form SINDy (WSINDy, see \cite{MessengerBortz2021JComputPhys,MessengerBortz2021SIAMMultiscaleModelSimul,messenger2021learning}), where integration against test functions replaces numerical differentiation, is a powerful method that is significantly more robust to noisy data, particularly in the context of PDEs. Furthermore, WSINDy's efficient convolutional formulation makes it a viable method for identifying PDEs under the constraints of limited memory capacity and computing power that exist in the online setting\footnote{The method developed here could also be adapted to the standard SINDy algorithm, however we choose to focus on the weak form for its demonstrated abilities to handle noisy data with low computational overhead.}.

The development of online algorithms is a relatively recent pursuit (\cite{zinkevich2003online,hazan2006efficient}), yet much progress has been made in applications to finance (\cite{hazan2009stochastic}), data processing (\cite{dixit2019online}), and predictive control (\cite{koller2018learning}) (see \cite{dall2020optimization,hoi_online_2021} for a recent surveys). In the context of sparse regression, several works have addressed online $\ell_1$-minimization and other methods of regularizing the $\ell_0$ pseudo-norm, although not in the context of learning dynamical systems (\cite{yang_stochastic_2020,zhai_tracking_2019,jialei_wang_online_2014,yuantao_gu_l_0_2009,kopsinis_online_2011,yilun_chen_sparse_2009}). To the best of our knowledge, neither SINDy nor WSINDy have been merged with an online learning algorithm for PDEs\footnote{There has, however, been work related to leveraging the equation learning ability of SINDy with Model Predictive Control (\cite{kaiser2018sparse}).}.

A successful approach for identifying PDEs and tracking parameters ``on the fly'' using multidimensional snapshots of data arriving sequentially over time would greatly benefit many areas of science and engineering. Possible paradigms in this online setting include identifying time-varying coefficients, SID in higher dimensions (where memory constraints require data to be streamed even for offline problems), and detecting changes in the dominant balance physics of the system, as terms become active or inactive dynamically. In this way, online sparse equation discovery has the potential to open doors to new application areas, and even improve performance of existing batch methods. 

We confront some of these challenges in this work by considering spatiotemporal dynamical systems and incoming data snapshots at every timestep. In the spirit of classical online algorithms, we develop an online WSINDy framework to this setting of streaming data with memory constraints by replacing full-data availability and batch optimization capabilities with data bursts and light-weight proximal gradient descent iterations to approximately solve the sparse regression problem. At each iteration we process only the incoming snapshot in time, and we do not assume the ability to compute least-squares projections apart from the initial guess. We focus on three prototypical systems, (1) the Kuramoto-Sivashinsky (KS) equation, which exhibits spatiotemporal chaos and thus has time-fluctuating Fourier content, (2) the nonlinear wave equation in a time-variable medium in two spatial dimensions, and (3) the linear wave equation in three spatial dimensions, a preliminary example of a system in higher dimensions.

\subsection{Notation}

Vector-valued objects will be bold and lower-case, $\xbf\in \Rbb^d$ for $d>1$, while multi-dimensional arrays will be bold and upper-case, $\Xbf\in \Rbb^{n_1\times \cdots \times n_d}$ for $n_i\in \Nbb$, $1\leq i\leq d$. To disambiguate between iteration and exponentiation, we refer to the $q$th element in a list of multi-dimensional arrays using superscripts in parentheses (e.g.\ $\xbf^{(q)}$ or $\Xbf^{(q)}$), whereas raising to the power $q$ (where applicable) is simply denoted $\Xbf^q$. Reference to an element within a multi-dimensional array is given as a subscript (e.g.\ $\xbf_i$ or $\Xbf_{i_1,\dots,i_d}$). For a matrix $\Gbf\in \Cbb^{m\times n}$, we denote by $\Gbf_{S}$ the restriction of $\Gbf$ to the columns in $S\subset \{1,\dots,n\}$. By some abuse of notation, $\Gbf_{S}^T = (\Gbf_{S})^T$. Similary, for a vector $\wbf\in \Cbb^n$, we let $\wbf_S \in \Rbb^{|S|}$ be the restriction of $\wbf$ to the entries in $S$, where $|S|$ denotes the number of elements of $S$. The complement of $S$ within $\{1,\dots,n\}$ is denoted $S^c$. All scalar-valued objects will be in lower-case, with iteration, set membership, etc. denoted by subscripts (i.e.\ $u_q$ is the $q$th element in the list $\{u_1,\dots, u_{q-1},u_q,u_{q+1},\dots\}$). 

\section{Problem Formulation}

We consider PDEs of the form
\begin{equation}\label{pdemodel}
D^{\pmb{\alpha}^{(0)}}u(\xbf,t) = \sum_{i,j=1}^{I,J} \wstar_{(i-1)J+j}(t) D^{\pmb{\alpha}^{(i)}}f_j(u(\xbf,t),\xbf), \quad (\xbf,t)\in \Omega\times[0,\infty),
\end{equation}
where $\Omega\subset\Rbb^d$ is a bounded open set. The operators $D^{\pmb{\alpha}^{(i)}}$ for $1\leq i \leq I$ represent any linear differential operator in the variables $(\xbf,t)\in \Rbb^{d+1}$, where $\pmb{\alpha}^{(i)}=(\pmb{\alpha}^{(i)}_1,\dots,\pmb{\alpha}^{(i)}_{d+1})$ is a multi-index such that 
\[D^{\pmb{\alpha}^{(i)}}v = \frac{\partial^{\pmb{\alpha}^{(i)}_1+\cdots+\pmb{\alpha}^{(i)}_d+\pmb{\alpha}^{(i)}_{d+1}}}{\partial \xbf_1^{\pmb{\alpha}^{(i)}_1}\cdots\partial \xbf_d^{\pmb{\alpha}^{(i)}_d}\partial t^{\pmb{\alpha}^{(i)}_{d+1}}}v.\]
In this work we consider left-hand side operators $D^{\pmb{\alpha}^{(0)}}$ to be either $\partial_t$ or $\partial_{tt}$, which are given in two spatial dimensions ($d=2$) by the multi-indices $\pmb{\alpha}^{(0)}=(0,0,1)$ and $\pmb{\alpha}^{(0)}=(0,0,2)$, respectively. The functions $f_j:\Rbb\times\Rbb^d \to \Rbb$, $1\leq j\leq J$, include all possible nonlinearities present in the model, and together with the linear operators $D^{\pmb{\alpha}^{(i)}}$ comprise the feature library $\Theta:=\{D^{\pmb{\alpha}^{(i)}}f_j\}_{i,j=1}^{I,J}$. The weight vector $\wstar(t) \in \Rbb^{IJ}$ is assumed to be {\it sparse} in $\Theta$ at each time $t$, and is allowed to vary in $t$.

We assume that at each time $t=k\Delta t$ for $k\in \Nbb$ and fixed timestep $\Delta t$ we are given a solution snapshot $\Ubf^{(t)}\in \Rbb^{n_1\times \cdots \times n_d}$ of the form
\begin{equation}\label{datassump}
\Ubf^{(t)} = u(\Xbf,t)+\ep
\end{equation}
where $u$ solves \eqref{pdemodel} for some weight vector $\wstar$ and $\Xbf \in \Rbb^{n_1\times \cdots \times n_d}$ is a fixed known spatial grid of points in $\Omega$ having $n_i$ points in the $i$th dimension and equal spacing $\Delta x$ in each dimension. Here $\ep$ represents i.i.d.\ mean-zero noise with fixed finite variance $\sigma^2$ associated with sampling the underlying solution $u(\xbf,t)$ at any point $\xbf\in \Omega$. We write $\Ubf = (\Ubf^{(0)},\Ubf^{(\Delta t)}, \dots, \Ubf^{(k\Delta t)},\dots)$ to denote the entire dataset in time. The problem is stated as follows.\\

\fbox{\begin{minipage}{0.92\textwidth}
\textbf{Problem:} Assume that a total of $K_\text{mem}$ snapshots $\{\Ubf^{(t-(K_\text{mem}-1)\Delta t)},\dots,\Ubf^{(t)}\}$ can be stored in memory at each time $t$ and that at time $t+\Delta t$ a new snapshot $\Ubf^{(t+\Delta t)}$ arrives, replacing the oldest snapshot in memory. Given the sampling model \eqref{datassump} for unknown $\sigma^2$, unknown ground truth PDE \eqref{pdemodel}, and fixed library $\Theta:=\{D^{\pmb{\alpha}^{(i)}}f_j\}_{i,j=1}^{I,J}$, solve for coefficients $\what^{(t)}$ such that $\sup_{t>0} \nrm{\what^{(t)}-\wstar(t)}$ is bounded.
\end{minipage}}

\section{Batch WSINDy}

In the batch setting, assuming $\wstar$ is constant in time, the weak-form sparse identification of nonlinear dynamics algorithm (WSINDy) proposed in \cite{MessengerBortz2021JComputPhys,MessengerBortz2021SIAMMultiscaleModelSimul} solves this problem efficiently by discretizing equation \eqref{pdemodel}, rewritten in its {\it convolutional weak form}:
\begin{equation}\label{weakform}
D^{\pmb{\alpha}^{(0)}}\psi*u(\xbf,t) = \sum_{i,j=1}^{I,J}\wstar_{(i-1)J+j} D^{\pmb{\alpha}^{(i)}}\psi*f_j(u,\cdot)(\xbf,t),
\end{equation}
where $\psi(\xbf,t)$ is any smooth function compactly supported in $\Omega\times [0,T]$ and convolutions are performed over space and time. For efficiency, the test function $\psi$ is chosen to be separable,
\begin{equation}\label{septest}
\psi(\xbf,t)=\phi_1(\xbf_1)\cdots\phi_d(\xbf_d)\phi_{d+1}(t) .
\end{equation}
For example, it can be chosen using the Fourier spectrum of the noisy data to mitigate high-frequency noise (see \cite{MessengerBortz2021JComputPhys}). Once $\psi$ is chosen, we discretize the problem by selecting a finite set of {\it query points} $\CalQ:=\{(\xbf^{(q)},t_q)\}_{q=1}^Q\subset \Omega\times (0,T)$ and evaluating \eqref{weakform} at $\CalQ$, replacing $u$ with the full dataset $\Ubf$. Convolutions can be efficiently computed using the fast Fourier transform (FFT), which, due to the compact support of $\psi$, is equivalent to the trapezoidal rule and is highly accurate in the noise-free case ($\sigma^2=0$). This gives us the linear system 
\[\bbf \approx \Gbf\wstar,\]
 where the $q$th entry of $\bbf$ is $\bbf_q = D^{\pmb{\alpha}^{(0)}}\psi*\Ubf(\xbf^{(q)},t_q)$ and $q$th entry of the $((i-1)J+j)$th column of $\Gbf$ is $\Gbf_{q,(i-1)J+j}=D^{\pmb{\alpha}^{(i)}}\psi*f_j(\Ubf,\cdot)(\xbf^{(q)},t_q)$. Using the assumption that $\wstar$ is sparse, we solve this linear system for $\what \approx \wstar$ by solving the sparse recovery problem
 \begin{equation}\label{sparserec}
\min_{\wbf\in \Rbb^{IJ}} F(\wbf;\lambda) = \min_{\wbf\in \Rbb^{IJ}} \frac{1}{2}\nrm{\Gbf\wbf-\bbf}_2^2 + \frac{1}{2}\lambda^2\nrm{\wbf}_0.
\end{equation}
The sparsity threshold $\lambda>0$ must be set by the user and is designed to strike a balance between fitting the data, associated with low residual $\nrm{\Gbf \wbf-\bbf}_2$, and finding a parsimonious model, indicated by low $\nrm{\wbf}_0$ (and its value is typically calibrated via cross-validation)  \cite{hastie2009elements,10.5555/2526243}.

In most cases, the columns of $\Gbf$ are highly correlated since they are each constructed from the same dataset $\Ubf$. This leads to many popular algorithms for solving \eqref{sparserec} performing poorly, such as convex relaxation using the $\ell_1$ norm or greedy search methods. In the batch setting, the following approach has proved to be successful under various noise levels and systems of interest. For $\lambda>0$ define the inner sequential thresholding step
\begin{equation}\label{MSTLS1}
\text{MSTLS}(\Gbf,\bbf; \lambda\,)\qquad \begin{dcases} \hspace{0.3cm}\wbf^{(0)} = \Gbf^\dagger \bbf \\ \hspace{0.43cm}\CalI^{(\ell)} = \{1\leq k\leq IJ\ :\ L_k(\lambda)\leq|\wbf^{(\ell)}_k|\leq U_k(\lambda)\} \\
\wbf^{(\ell+1)} = \argmin_{\supp{\wbf}\subset \CalI^{(\ell)}} \nrm{ \Gbf  \wbf-\bbf}_2^2.\end{dcases}
\end{equation}
Letting $\Gbf_{k}$ be the $k$th column of $\Gbf$, the lower and upper bounds are defined
\begin{equation}\label{MSTLSbnds} 
\begin{dcases} L_k(\lambda) =  \lambda\max\left\{1,\ \frac{\nrm{\bbf}}{\nrm{\Gbf_{k}}}\right\}\\
U_k(\lambda) =  \frac{1}{\lambda}\min\left\{1,\ \frac{\nrm{\bbf}}{\nrm{\Gbf_{k}}}\right\}\end{dcases}, \qquad 1\leq k\leq IJ.
\end{equation}
The sparsity threshold $\widehat{\lambda}$ is then selected as the smallest minimizer of the cost function \begin{equation}\label{lossfcn}
\CalL(\lambda) = \frac{\nrm{\Gbf(\wbf(\lambda)-\wbf(0))}_2}{\nrm{\Gbf\wbf(0)}_2}+\frac{\nrm{\wbf(\lambda)}_0}{IJ}
\end{equation}
where $\wbf(\lambda):=\text{MSTLS}(\Gbf,\bbf; \lambda\,)$. We find $\widehat{\lambda}$ via grid search and set $\what=\text{MSTLS}(\Gbf,\bbf;\widehat{\lambda})$ as the output of the algorithm. In words, this is a modified sequential thresholding algorithm with non-uniform thresholds \eqref{MSTLSbnds} chosen based on the norms of the underlying library terms $\Gbf_{(i-1)J+j}\approx D^{\pmb{\alpha}^{(i)}}\psi*f_j(u)$ relative to the response vector $\bbf \approx D^{\pmb{\alpha}^{(0)}}\psi*u$. The purpose of this is to (a) incorporate relative sizes of library terms $\Gbf_{k}\wstar_k$ along with absolute sizes of coefficients $\wstar$ in the thresholding step, and (b) choose $\lambda$ automatically. 

\section{Online WSINDy}

The online setting is defined by data snapshots $\Ubf^{(t)}$ arriving sequentially over time. An estimate $\what^{(t)}$ of the true parameters $\wstar(t)$ must be computed before the arrival of the next snapshot $\Ubf^{(t+\Delta t)}$ using only a fixed number $K_\text{mem}$ of previous snapshots. Without access to the full time series $\Ubf$, combined effects of the sample rate $\Delta t$, the number of snapshots $K_\text{mem}$, and the intrinsic timescales of the data determine the identifiability of the system: $\Delta t$ must be small enough to accurately compute time integrals, but large enough that the data $\Ubf$ is sufficiently dynamic over the time window $K_\text{mem} \Delta t$. Corruptions from noise have a greater impact because variance is not reduced by considering many samples in time, as was the case in the batch setting. Moreover, in realistic settings, solving for $\widehat{\wbf}^{(t)}$ before arrival of the next snapshot $\Ubf^{(t+\Delta t)}$ fundamentally limits the size of $(\Gbf,\bbf)$ and the number of iterations one may perform using any sparse solver.

The online setting is inherently restrictive, yet it appears well-suited for an important set of problems that are challenging offline and for settings where $\widehat{\wbf}^{(t)}$ must be obtained without revisiting past data. In particular, when the coefficient vector $\wstar$ varies over time, the library $\Theta$ must include time-dependent terms and may grow too large to successfully solve for an accurate sparse solution. Another issue arises with high-dimensional datasets (as in cosmology, turbulence, molecular dynamics, etc.), which cannot easily be processed in a single batch. In these cases an online approach is natural and advantageous even if solutions $\widehat{\wbf}^{(t)}$ are not themselves required ``online''.

For the online approach, at each time $t$ we seek to minimize the online cost function
\begin{equation}\label{sparserecOL}
\min_{\wbf\in \Rbb^{IJ}} F_t(\wbf;\lambda_t) = \min_{\wbf\in \Rbb^{IJ}} \frac{1}{2}\nrm{\Gbf^{(t)}\wbf-\bbf^{(t)}}_2^2 + \frac{1}{2}\lambda_t^2\nrm{\wbf}_0,
\end{equation}
where $(\Gbf^{(t)},\bbf^{(t)})$ is the linear system created from the $K_\text{mem}$ slices $\{\Ubf^{(t-(K_\text{mem}-1)\Delta t)},\dots,\Ubf^{(t)}\}$ at time $t$. Notice also that we allow $\lambda_t$ to change, as the initial guess $\lambda_0$ may not be optimal. In this online setting, we assume that we do not have the luxury of computing least-squares solutions (other than the initial guess), so we cannot use the approach outlined in \eqref{MSTLS1}-\eqref{lossfcn}, where \eqref{MSTLS1} requires multiple least-squares solves, and performing a grid search over $\lambda$ values requires multiple solves of \eqref{MSTLS1}. Hence, we consider the following online algorithm, which is simply the online proximal gradient descent combined with a decision tree update for $\lambda_t$ at each step:
\begin{equation}\label{WSINDyOL}
\begin{dcases}
\zbf^{(t)} = \widehat{\wbf}^{(t)} - \alpha_t(\Gbf^{(t)})^T\left(\Gbf^{(t)}\widehat{\wbf}^{(t)}-\bbf^{(t)}\right)\\
\what^{(t+\Delta t)} = H_{\lambda_t}\left(\zbf^{(t)}\right)\\
\lambda_{t+\Delta t}=\CalT(\lambda_t, \what^{(t+\Delta t)},\Delta \lambda,\lambda_{\max}).
\end{dcases}
\end{equation}
The hard thresholding operator $H_\lambda(\wbf)$ is the proximal operator of $\frac{1}{2}\lambda^2\nrm{\wbf}_0$ and is defined as
\begin{equation}
(H_\lambda(\wbf))_k = \begin{dcases} \wbf_k, &|\wbf_k|\geq \lambda \\0, & \text{otherwise}. \end{dcases}
\end{equation}
The map $\CalT$ updates $\lambda$ according to
\begin{equation}\label{lambdaupdate}
\CalT(\lambda_t, \what^{(t+\Delta t)}, \Delta \lambda,\lambda_{\max}) = \begin{dcases} (1-\Delta \lambda)\lambda_t, &\hspace{0.3cm} F_t(\what^{(t+\Delta t)},\lambda_t)>F_{t-1}(\widehat{\wbf}^{(t)},\lambda_t)\ \& \ S_{t+\Delta t}\subsetneq S_t \\ (1-\Delta \lambda )\lambda_t +\lambda_{\max}\Delta \lambda , &\begin{dcases} F_t(\what^{(t+\Delta t)},\lambda_t)>F_{t-1}(\widehat{\wbf}^{(t)},\lambda_t)\ \& \ S_{t}\subsetneq S_{t+\Delta t}\\
 F_t(\what^{(t+\Delta t)},\lambda_t)\leq F_{t-1}(\widehat{\wbf}^{(t)},\lambda_t)\ \& \  S_{t}=S_{t+\Delta t}.
\end{dcases} \\ \lambda_t, & \text{otherwise}.\end{dcases}
\end{equation}
In words, there are two possible updates to $\lambda_t$: a convex combination between $\lambda_t$ and $0$ and a convex combination between $\lambda_t$ and $\lambda_{\max}$. The former decreases $\lambda_t$ and occurs when library terms are thresholded to zero {\it and} the objective function $F_t$ increases. The latter increase $\lambda_t$ and occurs when either (a) library terms are added {\it and} $F_t$ increases or (b) the support set $S_t:=\supp{\widehat{\wbf}^{(t)}}$ doesn't change {\it and} $F_t$ does not increase.
At each step we set $\alpha_t = 1/\nrm{(\Gbf^{(t)})^T\Gbf^{(t)}_{S_t}}_2$, the optimal stepsize for pure gradient descent given the support $S_t$. As an initial guess we set $\what^{(0)}=\left(\Gbf^{(0)}\right)^\dagger\bbf^{(0)}$, which is the only least squares solve performed.

\begin{rmrk}
The value for $\lambda_t$ (and similarly for $\alpha_t$) can easily be replaced by a constant when additional knowledge is available (e.g. when $\wstar$ is known to satisfy certain bounds). While it is not common to update $\lambda_t$ over the course of an algorithm for sparse regression, it is well-known that picking $\lambda$ is problem specific and prone to errors particularly in the presence of noise (see \cite{MessengerBortz2021JComputPhys} for a discussion). The update policy given by $\CalT$ encodes simple objectives of any minimization algorithm for \eqref{sparserecOL} and works in all examples presented. A full investigation is a topic for future work.
\end{rmrk}

\begin{rmrk} 
Similar to the batch case, we find that non-uniform thresholding greatly improves results. For brevity, we include in Appendix \ref{App:tech} a description of how non-uniform thresholds such as \eqref{MSTLSbnds} can be incorporated in the online framework. We also note that the theoretical results in the next section carry over analogously in the non-uniform thresholding case.  \end{rmrk}

\subsection{Regret and Fixed Point Analysis}

The behavior of the online algorithm is in large part dictated by the behavior of the batch proximal gradient descent method. To the best of our knowledge the proximal gradient descent algorithm applied to the $\ell_0$ norm has not been well studied. First we review useful properties relating solutions of \eqref{sparserec} and stationary points of the proximal gradient descent algorithm \eqref{WSINDyOL} in the offline case and for fixed $\lambda$. We then use these results to bound the dynamic online regret, which we define as 
\begin{equation}\label{regret}
Reg_D(T):= \sum_{t=0}^T F_t(\widehat{\wbf}^{(t)};\lambda_t)-F_t(\wstar(t);\lambda_t),
\end{equation}
where $\wstar(t)$ is a global minimizer of $F_t(\wbf,\lambda_t)$. In particular, we first have the following:
\begin{lemm}\label{lemm}
Consider $\wbf$ such that one of the following holds:
\begin{enumerate}[label=(\roman*)]
\item $\wbf$ is a local minimizer of \eqref{sparserec}
\item $\wbf = H_\lambda\left(\wbf - \Gbf^T\left(\Gbf \wbf - \bbf\right)\right)$
\item With $S=\supp{\wbf}$, we have that $\wbf_S \in \argmin_\zbf\nrm{\Gbf_{S}\zbf-\bbf}_2^2$ and
\[\max_{i\in S^c}\left\vert \Gbf^T_{i}(\Gbf\wbf-\bbf) \right\vert<\lambda \leq \min_{i\in S} \left\vert \wbf_i\right\vert.\]
\end{enumerate}
Then it holds that $(ii)\iff(iii) \implies (i)$. Moreover, if $\wbf$ a global minimizer, then $(i)\implies (iii)$.
\end{lemm}

The proof of Lemma \ref{lemm} can be found in Appendix \ref{app:lemm}. Lemma \ref{lemm} implies in particular that fixed points of the batch proximal gradient descent algorithm are local minimizers of $F(\wbf;\lambda)$, and moreover that fixed points satisfy a necessary condition for global optimality given by $(iii)$. If in addition a fixed point $\wbf$ with $\supp{\wbf}=S$ satisfies that $\Gbf_{S}^T\Gbf_{S}$ is full rank,  then $\wbf_S = \Gbf_{S}^\dagger\bbf$ is the unique least squares solution over the columns in $S$. In \cite{nikolova2013description} it is shown that this is sufficient for $\wbf$ to be a {\it strict} local minimizer, which then implies that $\wbf$ is the {\it unique} global minimizer {\it restricted to the set $S$}. Also in \cite{nikolova2013description} is an extensive treatment of global minimizers of $F(\wbf;\lambda)$, where it is shown that apart from a measure-zero set of linear systems $(\Gbf,\bbf)$, the global minimizer is unique. We use this to bound the dynamic regret below.

\begin{thm}
Let $\sigma_{1,t}$ and $\sigma_{n,t}$ denote the first and last singular values of the matrix $\Gbf^{(t)} \in \Rbb^{m\times n}$. Assume the following: $\max_t\lambda_t\leq \overline{\lambda}<\infty$, $\min_t \sigma_{n,t} \geq \overline{\sigma}_{\min}>0$, $\max_t \sigma_{1,t} \leq \overline{\sigma}_{\max}$, and $\alpha_t<\sigma_{\max}^{-2}$. In addition, assume that the global minimizer $\wstar(t)$ of $F_t(\wbf;\lambda_t)$ is unique for every $t$ and satisfies $|S^\star_t|\geq \overline{s}>0$ where $S^\star=\supp{\wstar(t)} $. Finally, assume that the tracking gap is globally bounded: $\nrm{\wstar(t)-\wstar(t+1)}_2:=d_t\leq \overline{d}$. Then the dynamic regret \eqref{regret} grows at-worst linearly:
\[Reg_D(T) \leq C_1+C_2T\]
for some $C_1>0$ and $C_2>0$. In particular, $\frac{1}{T}Reg_D(T)$ remains bounded.
\end{thm}
The constants $C_1$ and $C_2$ are specified in the proof, which is presented next. 
\begin{proof}
First we decompose $F_t(\wbf;\lambda_t) = g_t(\wbf)+h_t(\wbf) = \nrm{\Gbf^{(t)}\wbf-\bbf^{(t)}}_2^2+\lambda^2\nrm{\wbf}_0$. We can bound the difference in $g_t$ as follows:
\begin{align*}
g_t(\wbf^{(t)})-g_t(\wstar(t)) &= \nrm{\Gbf^{(t)}\wbf^{(t)}}_2^2-\nrm{\Gbf^{(t)}\wstar(t)}_2^2-2\lan \bbf^{(t)},\Gbf^{(t)}\left(\wbf^{(t)}-\wstar(t)\right)\ran \\
&=\nrm{\Gbf^{(t)}(\wbf^{(t)}-\wstar(t))}_2^2-2\lan \Gbf^{(t)}(\wbf^{(t)}-\wstar(t)),\Gbf^{(t)}\wstar(t)-\bbf^{(t)}\ran\\
\intertext{taking $\left\vert\cdot\right\vert$ of both sides and noting from the Lemma that $\nrm{(\Gbf^{(t)})^T(\Gbf^{(t)}\wstar(t)-\bbf^{(t)})}_\infty< \lambda_t$ implies that}
g_t(\wbf^{(t)})-g_t(\wstar(t))&\leq \overline{\sigma}_{\max}^2\nrm{\wbf^{(t)}-\wstar(t)}_2^2+2\lambda_t\sqrt{|(S^\star)^c|}\nrm{\wbf^{(t)}-\wstar(t)}_2\\
&\leq \overline{\sigma}_{\max}^2\nrm{\wbf^{(t)}-\wstar(t)}_2^2+2\overline{\lambda}\sqrt{n-\overline{s}}\nrm{\wbf^{(t)}-\wstar(t)}_2.
\end{align*}
For $h_t$ we have simply
\[\left\vert h_t(\wbf^{(t)})-h_t(\wstar(t))\right\vert = \lambda_t^2\left\vert |S_t|-|S_t^\star|\right\vert \leq \overline{\lambda}^2(n-\overline{s}).\]
For any vectors $\xbf,\ybf\in \Rbb^n$, it holds that
\[\nrm{H_\lambda(\xbf)-H_\lambda(\ybf)}\leq \nrm{\xbf-\ybf}_2 + \lambda \sqrt{|S_x\triangle S_y|}\]
where $S_x\triangle S_y = (S_x\cap S_y^c)\cup (S_y\cap S_x^c)$ is the symmetric difference of the sets $S_x=\supp{H_\lambda(\xbf)}$ and $S_y = \supp{H_\lambda(\ybf)}$. This implies, together with stationarity of $\wstar(t)$,
\begin{align*}
&\nrm{\wbf^{(t+\Delta t)}-\wstar(t)}_2 \\
&=\nrm{H_{\alpha_t\lambda_t}\left(\wbf^{(t)}-\alpha_t(\Gbf^{(t)})^T\left(\Gbf^{(t)}\wbf^{(t)}-\bbf^{(t)}\right)\right)-H_{\alpha_t\lambda_t}\left(\wstar(t)-\alpha_t(\Gbf^{(t)})^T\left(\Gbf^{(t)}\wstar(t)-\bbf^{(t)}\right)\right)}_2\\
&\leq \nrm{\left(\Ibf-\alpha_t(\Gbf^{(t)})^T\Gbf^{(t)}\right)(\wbf^{(t)}-\wstar(t))}_2 + \alpha_t\lambda_t\sqrt{|S_{t+1}\triangle S_t^\star|}\\
&\leq \max\left(|1-\alpha_t\sigma_{1,t}^2|,|1-\alpha_t\sigma_{n,t}^2|\right)\nrm{\wbf^{(t)}-\wstar(t)}_2 + \alpha_t\lambda_t \sqrt{|S_{t+1}\triangle S_t^\star|}\\
&:=\rho_t\nrm{\wbf^{(t)}-\wstar(t)}_2 + \alpha_t\lambda_t \sqrt{|S_{t+1}\triangle S_t^\star|}.
\end{align*}
Using that  $\nrm{\wbf^{(t+\Delta t)}-\wstar(t+1)}_2\leq \nrm{\wbf^{(t+\Delta t)}-\wstar(t)}_2+d_t$, we have the recurrence relation 
\[\nrm{\wbf^{(t+\Delta t)}-\wstar(t+1)}_2 \leq \rho_t\nrm{\wbf^{(t)}-\wstar(t)}_2+d_t + \alpha_t\lambda_t \sqrt{|S_{t+1}\triangle S_t^\star|},\]
where, by assumptions on $\sigma_{1,t},\sigma_{n,t}$ and $\alpha_t$, it holds that $\max_t\rho_t\leq \overline{\rho}$ for some $\overline{\rho}<1$, hence we get the bound
\[\nrm{\wbf^{(t)}-\wstar(t)}_2\leq \overline{\rho}^t\nrm{\wbf^{(0)}-\wstar(0)}_2 + (\overline{d}+\overline{\alpha}\overline{\lambda} \sqrt{n})\sum_{s=0}^t\overline{\rho}^s \leq \overline{\rho}^t\nrm{\wbf^{(0)}-\wstar(0)}_2 +\frac{\overline{d}+\overline{\alpha}\overline{\lambda} \sqrt{n}}{1-\overline{\rho}}.\]
We note in passing that this implies a uniform error bound on $\nrm{\wbf^{(t)}-\wstar(t)}_2$ which asymptotically depends only on the tracking gap $d_t$ and the support difference $|S_t\triangle S_t^\star|$. Finally, using this bound and previous calculations for $g$ and $h$, we get 
\[Reg_D(T) \leq \tilde{C}_1\sum_{t=0}^T\overline{\rho}^t + C_2 T \leq \frac{\tilde{C}_1}{1-\overline{\rho}} +C_2 T\]
where 
\[\tilde{C}_1 = \overline{\sigma}_{\max}^2\overline{\rho}\nrm{\wbf^{(0)}-\wstar(0)}_2^2 + \frac{2\nrm{\wbf^{(0)}-\wstar(0)}_2}{1-\overline{\rho}}\left(\overline{\lambda}\sqrt{n-\overline{s}}+\frac{(\overline{d}+\overline{\alpha}\overline{\lambda} \sqrt{n})\overline{\sigma}^2_{\max}}{1-\overline{\rho}}\right)\] 
\[C_2 = \overline{\lambda}^2(n-\overline{s})+2\overline{\lambda}\sqrt{n-\overline{s}}\frac{(\overline{d}+\overline{\alpha}\overline{\lambda} \sqrt{n})}{1-\overline{\rho}}+\overline{\sigma}_{\max}^2\left(\frac{\overline{d}+\overline{\alpha}\overline{\lambda} \sqrt{n}}{1-\overline{\rho}}\right)^2\]
\[=(\overline{\sigma}_{\max}^2-1)\left(\frac{\overline{d}+\overline{\alpha}\overline{\lambda} \sqrt{n}}{1-\overline{\rho}}\right)^2 + \left(\frac{\overline{d}+\overline{\alpha}\overline{\lambda} \sqrt{n}}{1-\overline{\rho}}+\overline{\lambda}\sqrt{n-\overline{s}}\right)^2.\]
This completes the proof.
\end{proof}
The above result establishes that $Reg_D(T)$ increases at-worst asymptotically at the same rate as online gradient descent applied to the time-varying ordinary least squares problem (\cite{zinkevich2003online}), although with constant depending on both the tracking error $\overline{d}$ and a sparsity factor, which in our estimate takes the form $\overline{\lambda}\sqrt{n}$ ($n=IJ$ being the number of columns in $\Gbf^{(t)}$). Below we examine only $S^\star_t=S^\star$ fixed, and find that over a wide range of parameters the correct support is found in finitely many iterations, leading to scenarios where the dynamic regret depends only on $\overline{d}$ for large $T$ (see Figure \ref{fig:W2Dc} for a visualization of this case for the time-varying wave equation).

\section{Numerical Experiments}

Our primary focuses are the performance of the algorithm as a function of the number of snapshots $K_\text{mem}$ allowed in memory and the sensitivity of the algorithm to noise. We examine the following three examples which display a range of dynamics over one to three spatial dimensions: the Kuramoto-Sivashinsky equation in 1D, a time-varying nonlinear wave equation in 2D, and the linear wave equation in 3D. We abbreviate each by KS, W2D and W3D. For each experiment we simulate a noise-free solution $\Ubf_{exact}$ to the given PDE over a long time horizon. We then add i.i.d.\ Gaussian noise with mean zero and standard deviation $\sigma=\sigma_{NR}\nrm{\Ubf^\star}_{rms}$ to each data point for a range of {\it noise ratios}\footnote{Note that $\sigma_{NR}$ is approximately equal to the ratio $\nrm{\ep(:)}_2/\nrm{\Ubf_{exact}(:)}_2$ of the noise to the true data, where ``$\Ubf_{exact}(:)$'' denotes $\Ubf_{exact}$ stretched into a column vector.} $\sigma_{NR}$. After an offline phase where a least squares solution is found from the first $K_\text{mem}$ snapshots, we feed in one new snapshot at each time $t$ and apply the online algorithm \eqref{WSINDyOL}.

\subsubsection{Algorithm Hyperparameters}

We fix as many hyperparameters across examples as possible, and differences are summarized in Table \ref{stats}.  In all examples we fix the sparsity threshold update to $\Delta \lambda=0.1$, the initial sparsity threshold to $\lambda_0=0.0001$, and the maximum sparsity threshold to $\lambda_{\max}=0.1$. For the library we use 
\[\Theta = \{\partial_{\xbf_i}^k (u^j)\}, \qquad 1\leq i\leq d,\ 0\leq k\leq 4,\ 0\leq j\leq 4\] 
in other words all spatial derivatives up to degree 4 of monomials up to degree 4 of the data (excluding mixed derivatives). For direct comparison of the effects of $K_\text{mem}$ and $\sigma_{NR}$ across examples, we fix the test function $\psi$ in the representation \ref{septest} so that 
\begin{align}
\phi_i(\xbf_i) &= \left(1-\left(\frac{\xbf_i}{21\Delta x}\right)^2\right)^{11}_+, \qquad 1\leq i \leq d\\
\intertext{and}
\phi_{d+1}(t) &= \left(1-\left(\frac{t}{(K_\text{mem}-1)\Delta t/2}\right)^2\right)^{9}_+, \label{timephi}
\end{align}
where $(z)_+:=\max\{z,0\}$. In this way $\psi$ is supported on $2\times 21+1=43$ points in each spatial dimension and $K_\text{mem}$ points in time, although note that ($\Delta x, \Delta t$) change across examples. Since $\phi_{d+1}(t)$ is supported on $K_\text{mem}$ points, there is only one integration in time at each iteration, so that the query points are given by $\CalQ = \{(\xbf^{(q)},t_q)\}_{q=1}^Q = \CalQ_\xbf\times \{t-(K_\text{mem}-1)\Delta t/2\}$ where for each example $\CalQ_\xbf$ is fixed across all values of $K_\text{mem}$ and $\sigma_{NR}$. We take $\CalQ_\xbf\subset \Xbf$ to be equally-spaced and such that the linear system $(\Gbf^{(t)},\bbf^{(t)})$ contains less than 10,000 rows (see Table \ref{stats} for exact dimensions). Online iteration times are reported below for computations performed on a laptop with 1.7GHz base clockspeed AMD Ryzen 7 pro 4750u processor and 38.4 GB of RAM.

\begin{rmrk}
 
By defining the temporal test function $\phi_{d+1}(t)$ to depend on $K_\text{mem}$ according to \eqref{timephi}, the implied strategy is that increasing $K_\text{mem}$ (keeping more snapshots in memory) leads to more accurate integration in the time domain. One could instead fix the test function
\[\phi_{d+1}(t) = \left(1-\left(\frac{t}{m\Delta t}\right)^2\right)^9_+\] 
for some $m\leq (K_\text{mem}-1)/2$ for all $K_\text{mem}$ considered, leading to a fixed integration window of length $2m+1$ in time. Increasing $K_\text{mem}$ would then allow for more integrations in time (i.e.\ a larger set of query points $\CalQ$), adding rows to the linear system $(\Gbf^{(t)}, \bbf^{(t)})$. Our chosen strategy fixes the dimensions of $(\Gbf^{(t)}, \bbf^{(t)})$, leading to a more direct comparison across examples. We leave this trade-off between the number of time integrations and the accuracy of time integrations to future work.
\end{rmrk}

\begin{table}
\begin{center}
\begin{tabular}{|c|c|c|c|c|}
\hline    & dims$(\Xbf)$ & $T$ & dims$(\Gbf^{(t)})$ & $(\Delta x, \Delta t)$\\
\hline KS & $256 \times 1$ & $3946$ & $214 \times 21$ & $(0.939,0.586)$ \\
\hline W2D & $129\times 403$ & $1639$ & $7964 \times 37$ & $(0.0156,0.0122)$\\
\hline W3D & $128\times 128\times 128$ & $960$ & $8192 \times 53$ & $(0.0491,0.0122)$\\
\hline
\end{tabular}
\end{center}
\caption{Resolution and dimensions of datasets used in examples.}
\label{stats}
\end{table}

\subsubsection{Performance Analysis}

We are concerned with the ability of the algorithm to recover the support of the true model coefficients $S^\star:=\supp{\wstar}$ as well as the accuracy of $\widehat{\wbf}^{(t)}$ over time, depending primarily on the number $K_\text{mem}$ of solution snapshots allowed in memory and the noise level $\sigma_{NR}$ corrupting the data. To assess support recovery, we measure the {\it true positivity ratio} (TPR)
\[\text{TPR}(\widehat{\wbf}^{(t)}):=\frac{\text{TP}(\widehat{\wbf}^{(t)})}{\text{TP}(\widehat{\wbf}^{(t)})+\text{FP}(\widehat{\wbf}^{(t)})+\text{FN}(\widehat{\wbf}^{(t)})}\]
where $\text{TP}(\widehat{\wbf}^{(t)}):=|S_t\cap S^\star|$ is the number of correctly identified nonzero coefficients, $\text{FP}(\widehat{\wbf}^{(t)}):=|S_t\cap (S^\star)^c|$ is the number of falsely identified nonzero coefficients, and $\text{FN}(\widehat{\wbf}^{(t)}):=|S_t^c\cap S^\star |$ is the number of falsely identified zero coefficients. A TPR of 1 indicates successful support recovery, while $\text{TPR}=0.75$ indicates 3/4 terms were correctly identified, and so on. We measure the accuracy of $\widehat{\wbf}^{(t)}$ in the relative $\ell_2$-norm:
\[E_2(\widehat{\wbf}^{(t)}) := \nrm{\widehat{\wbf}^{(t)}-\wstar(t)}_{2}/\nrm{\wstar(t)}_{2}.\]
We report the results of $\text{TPR}(\widehat{\wbf}^{(t)})$ and $E_2(\widehat{\wbf}^{(t)})$ averaged over 100 instantiations of noise.


\subsection{Kuramoto-Sivashinsky (KS)}

\begin{equation}\label{KSeq}
\partial_t u = -\partial_x\left(u^2\right)-\partial_{xx} u-\partial_{xxxx}u.
\end{equation}

The Kuramoto-Sivashinsky (KS) equation is challenging because the solution exhibits spatiotemporal chaos and so has a Fourier spectrum that varies in time. This leads to potentially different dynamics at each timestep in the online learning perspective. The PDE also has a 4th-order derivative in space which is difficult to compute accurately and to identify via sparse regression, especially when noise is present. We simulate the solution using a high-order method (accurate to 6-7 digits) and use a dataset of $256\times 3496$ points in space and time at resolution $(\Delta x,\Delta t)= (0.393,0.586)$. Online iterations take less than  $0.01$ seconds, which includes building the linear system $(\Gbf^{(t)},\bbf^{(t)})$, which is the most costly step.

\begin{figure}
\centering
\begin{tabular}{cc}
\includegraphics[trim={0 0 20 0},clip,width=0.4\textwidth]{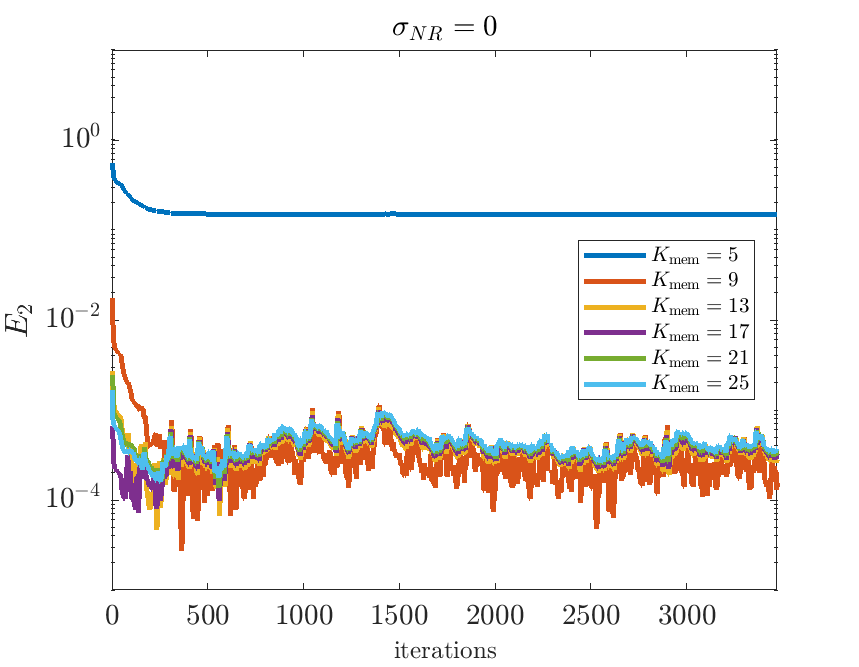} & 
\includegraphics[trim={0 0 20 0},clip,width=0.4\textwidth]{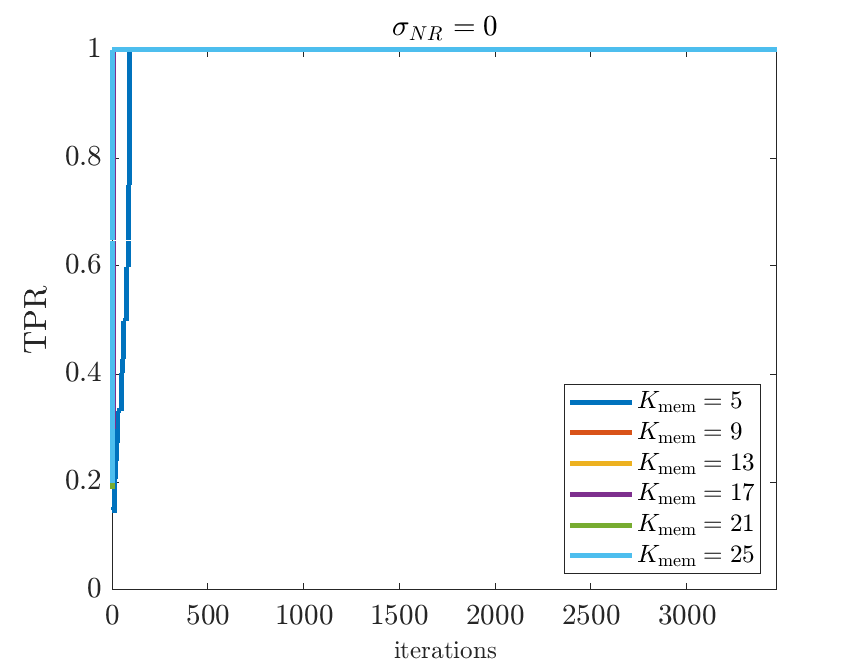} \\ 
\includegraphics[trim={0 0 20 0},clip,width=0.4\textwidth]{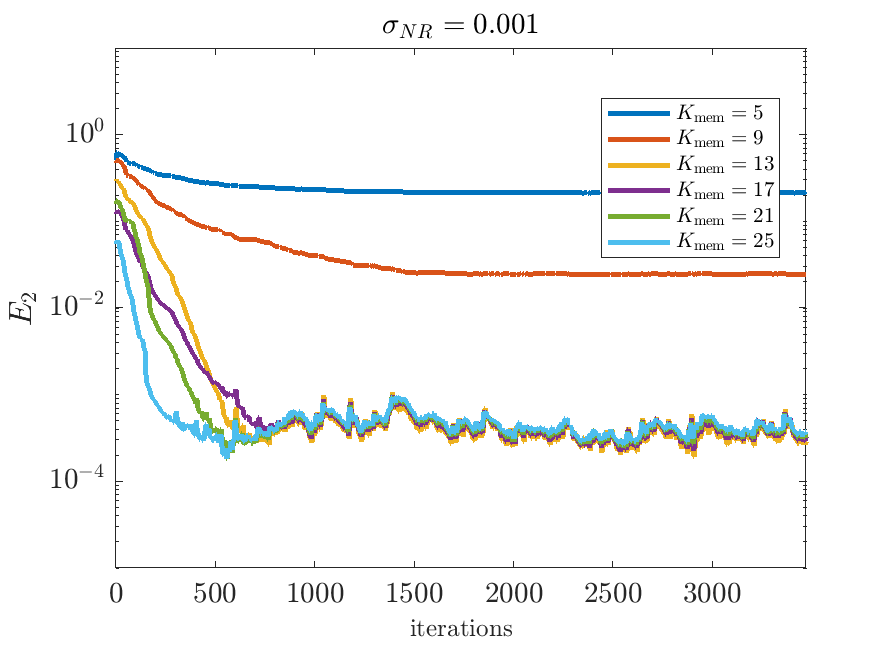} & 
\includegraphics[trim={0 0 20 0},clip,width=0.4\textwidth]{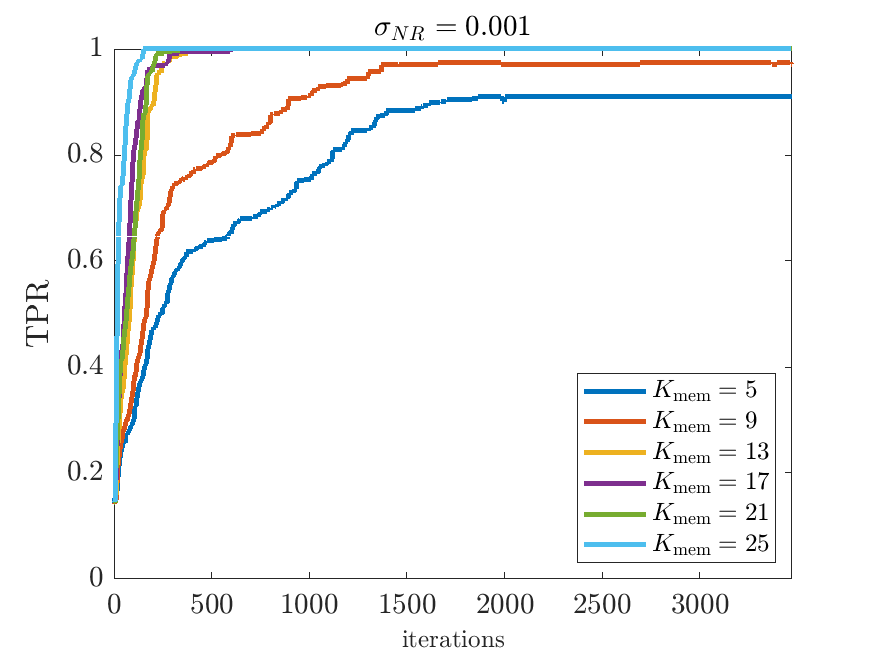} \\ 
\includegraphics[trim={0 0 20 0},clip,width=0.4\textwidth]{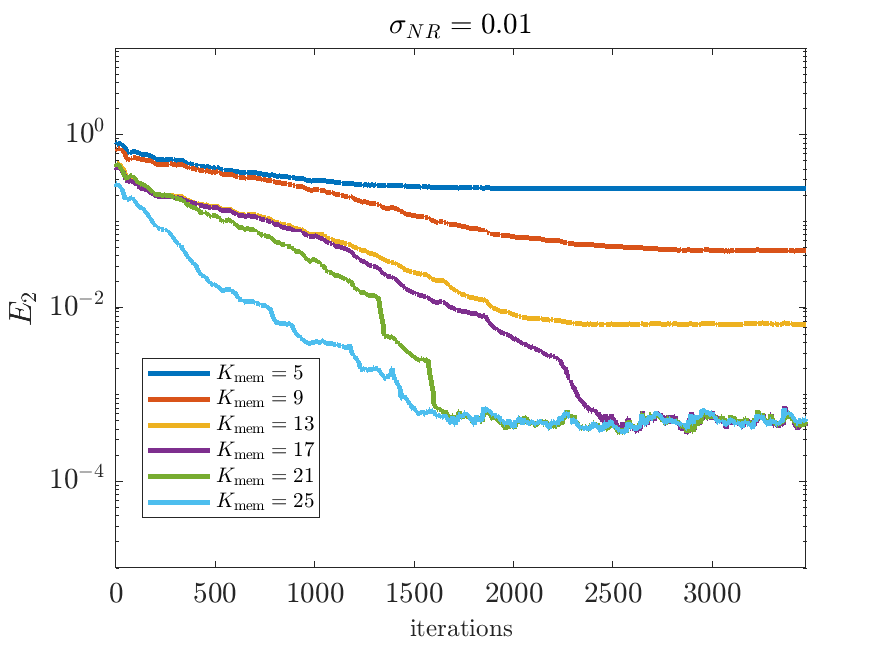} & 
\includegraphics[trim={0 0 20 0},clip,width=0.4\textwidth]{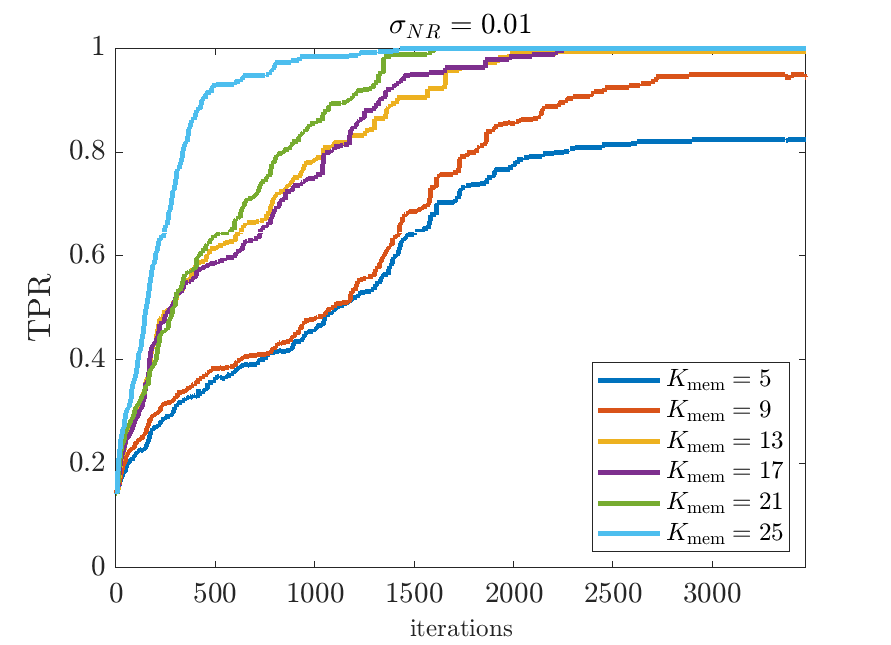} \\ 
\includegraphics[trim={0 0 20 0},clip,width=0.4\textwidth]{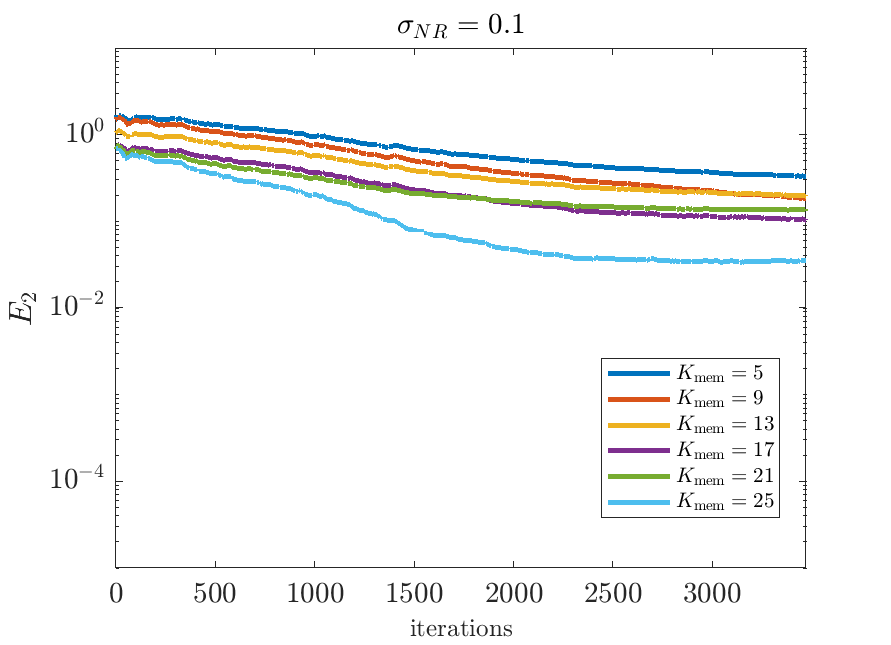} & 
\includegraphics[trim={0 0 20 0},clip,width=0.4\textwidth]{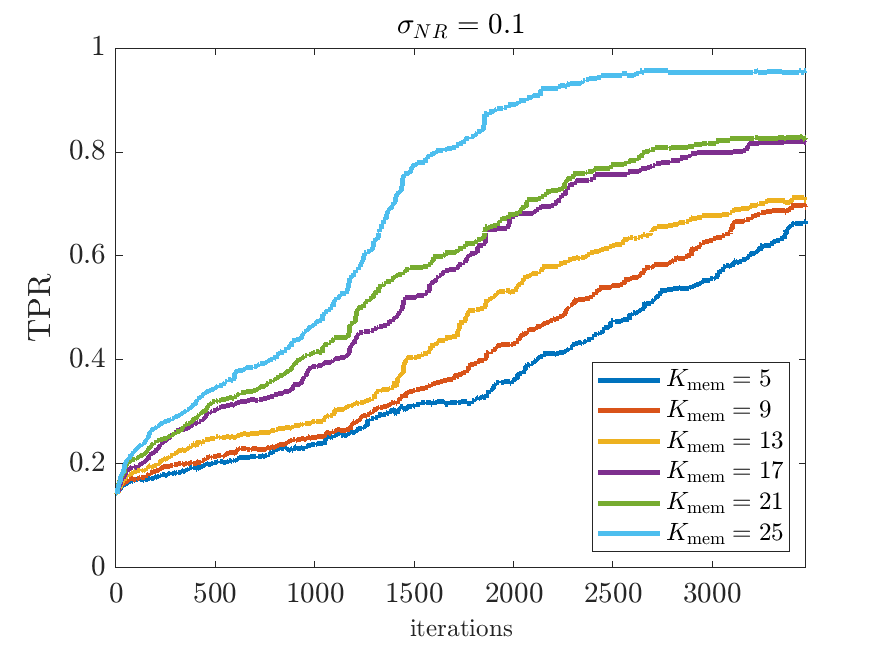}
\end{tabular}
\caption{Online identification of the Kuramoto-Sivashinsky equation \eqref{KSeq} for $K_\text{mem}\in\{5,9,13,17,21,25\}$ and (top to bottom) $\sigma_{NR}\in\{0,0.001,0.01,0.1\}$. Left: average coefficient error $E_2(\what^{(t)})$. Right: average total positivity ratio TPR$(\what^{(t)})$.}
\label{fig:KS}
\end{figure}

In Figure \ref{fig:KS} the average evolution of $E_2(\widehat{\wbf}^{(t)})$ and $\text{TPR}(\widehat{\wbf}^{(t)})$ is depicted for various noise levels $\sigma_{NR}$ and memory capacities $K_\text{mem}$. The system is correctly identified for all trials when $K_\text{mem}\in \{13,17,21,25\}$ and $\sigma_{NR}\in \{0,0.001,0.01\}$, with relative errors $E_2$ less than $10^{-2}$ once the system is identified. For larger noise $\sigma_{NR}=0.1$, results stagnate at sub-optimal values, indicating that more data is needed to identify the system (note that $\Gbf^{(t)}$ only has 214 rows). With $K_\text{mem}=5$ we recover the correct system only in the noiseless case ($\sigma_{NR}=0$), indicating that 5 points in time does not result in accurate resolution of the dynamics. 

\subsection{Variable-medium nonlinear wave equation in 2D (W2D)}

\begin{equation}\label{SGeq}
\partial_{tt} u = c(t)\left(\partial_{xx}u+\partial_{yy}u\right) -u^3
\end{equation}

\begin{figure}
\centering
\begin{tabular}{cc}
\footnotesize{$K_\text{mem}=17,\ \sigma_{NR}=0.01$} & \footnotesize{$K_\text{mem}=25,\ \sigma_{NR}=0.01$} \\
\includegraphics[trim={0 0 20 10},clip,width=0.4\textwidth]{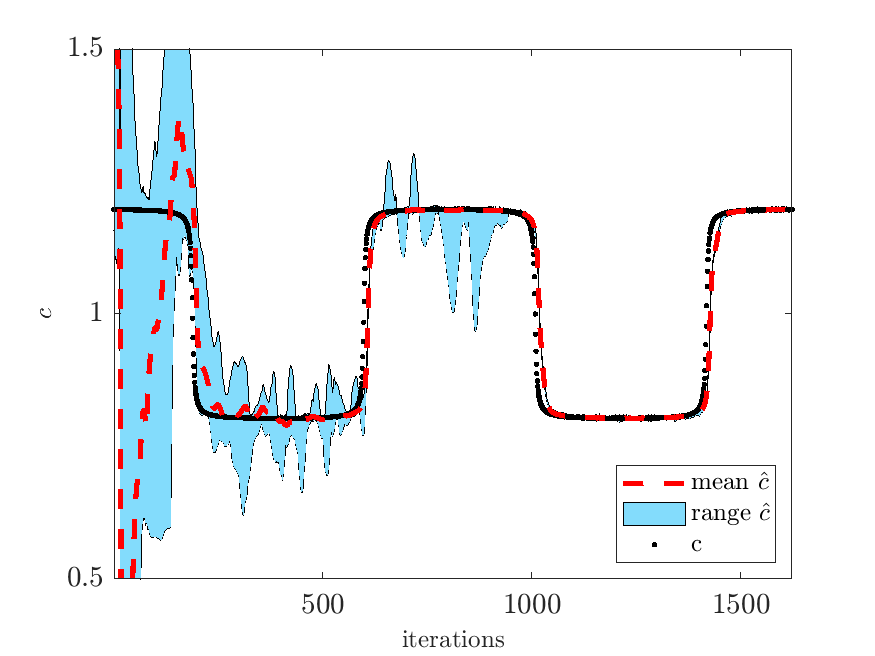} & 
\includegraphics[trim={0 0 20 10},clip,width=0.4\textwidth]{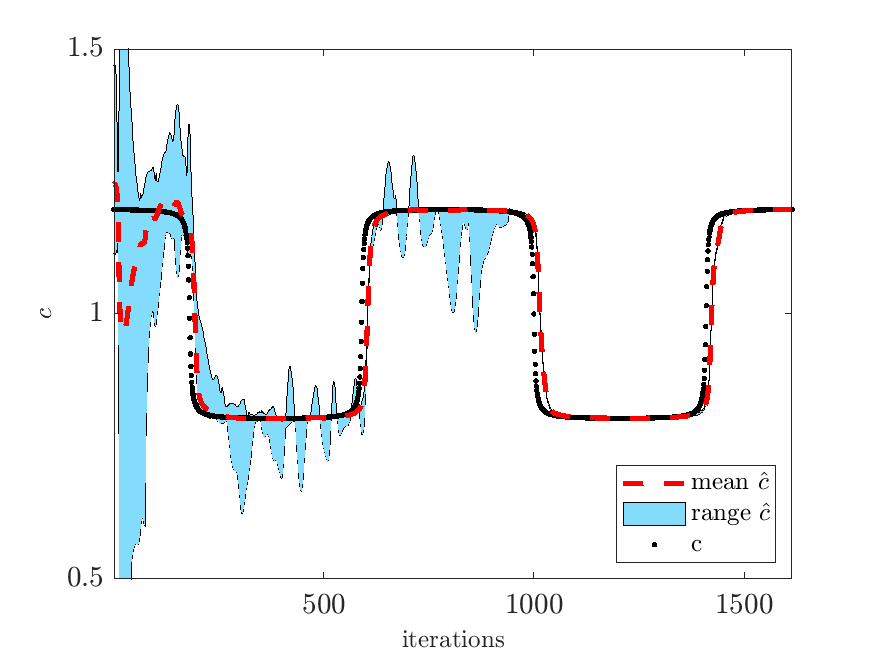} \\ 
\footnotesize{$K_\text{mem}=17,\ \sigma_{NR}=0.1$} & \footnotesize{$K_\text{mem}=25,\ \sigma_{NR}=0.1$} \\
\includegraphics[trim={0 0 20 10},clip,width=0.4\textwidth]{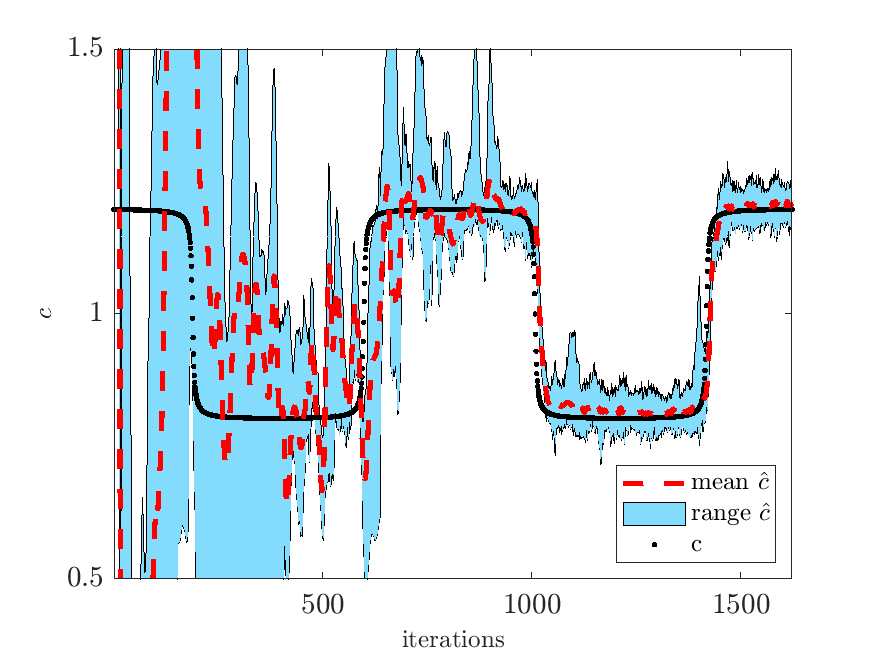} & 
\includegraphics[trim={0 0 20 10},clip,width=0.4\textwidth]{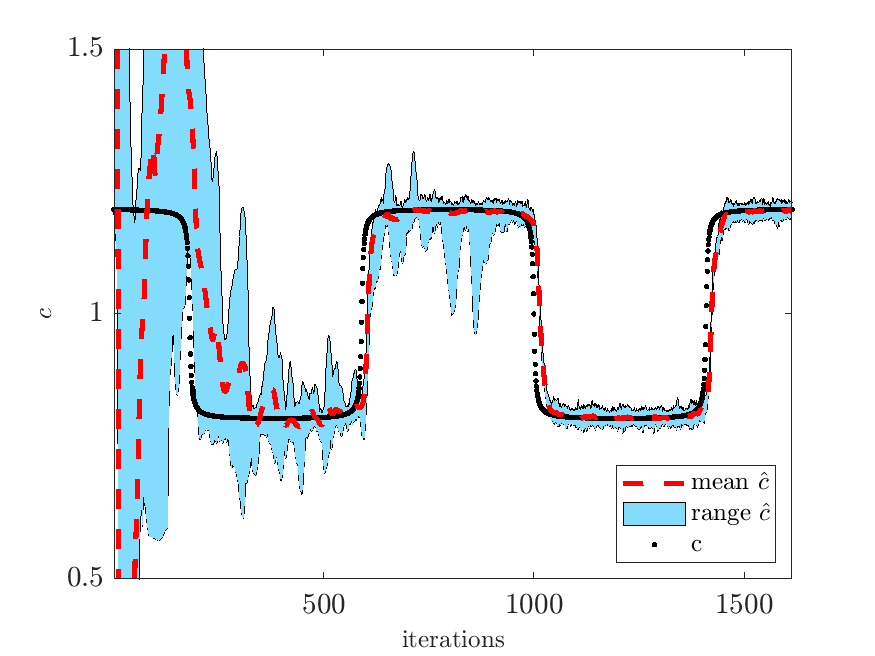}
\end{tabular}
\caption{Online estimation of the wavespeed $c(t)$ (shown in black) for PDE \eqref{SGeq}. The average learned wavespeed $\widehat{c}(t)$ is shown in red while the blue shaded region shows the maximum and minimum values attained over all 100 trials. Notice the accuracy for later iterations when $\sigma_{NR}=0.01$, and the reduction in variance moving from $K_\text{mem}=17$ to $K_\text{mem}=25$ when $\sigma_{NR}=0.1$.}
\label{fig:W2Dc}
\end{figure}

We examine a variable-medium nonlinear wave equation in 2D, given by equation \eqref{SGeq}, where the variable medium is modeled by the time-varying wavespeed
\[c(t) = 1 + (0.2)\frac{2}{\pi}\arctan(40\cos(2\pi(0.1)t)),\]
The wavespeed is a smoothed square wave and represents a system with abrupt speed modulation (see Figure \ref{fig:W2Dc} for depictions). We simulate the solution using a Fourier $\otimes$ Legendre spectral method in space with leap-frog timestepping. The exact data $\Ubf_{exact}$ has dimensions $129\times 403 \times 1639$ in $(x,y,t)$ with resolution $(\Delta x,\Delta t)= (0.0156,0.0122)$. Each snapshot $\Ubf^{(t)}$ is $0.42$ megabytes (Mb) and online iterations take approximately $0.08$ seconds.

\begin{figure}
\centering
\begin{tabular}{cc}
\includegraphics[trim={0 0 20 0},clip,width=0.4\textwidth]{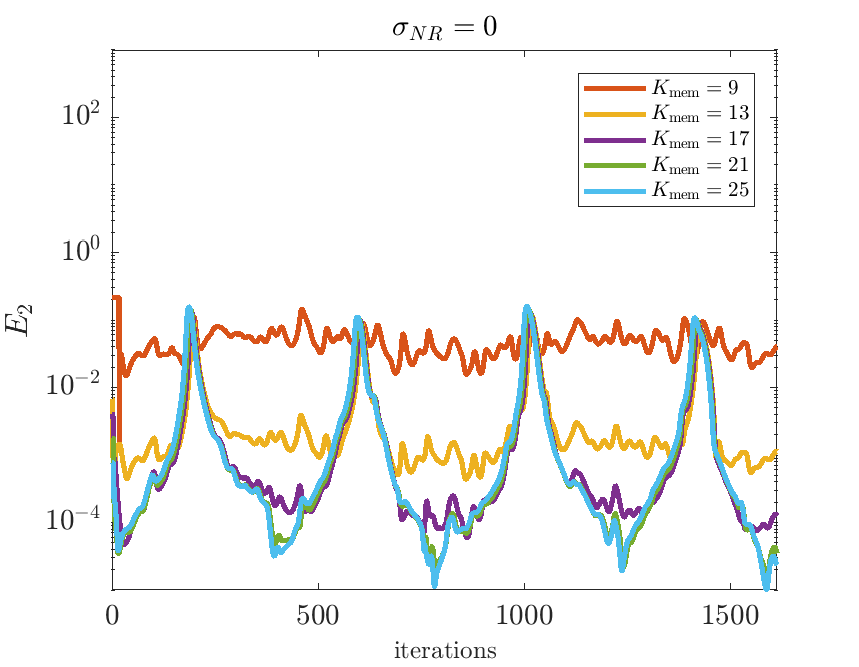} & 
\includegraphics[trim={0 0 20 0},clip,width=0.4\textwidth]{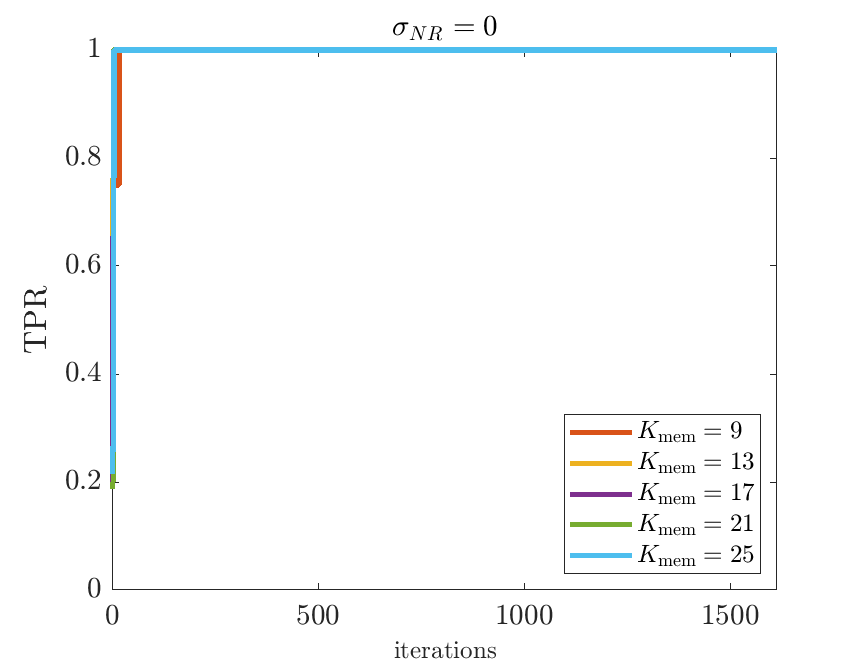} \\ 
\includegraphics[trim={0 0 20 0},clip,width=0.4\textwidth]{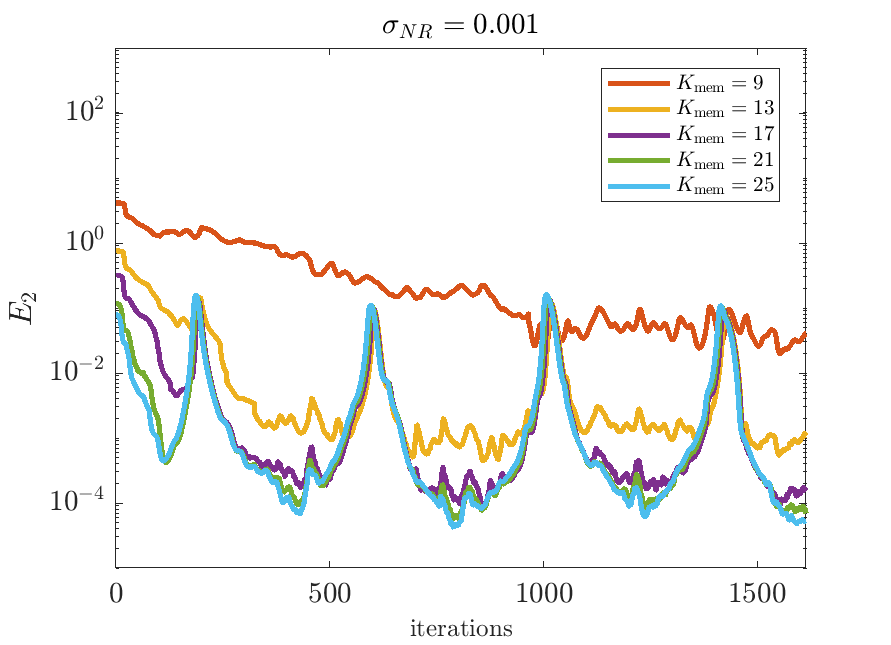} & 
\includegraphics[trim={0 0 20 0},clip,width=0.4\textwidth]{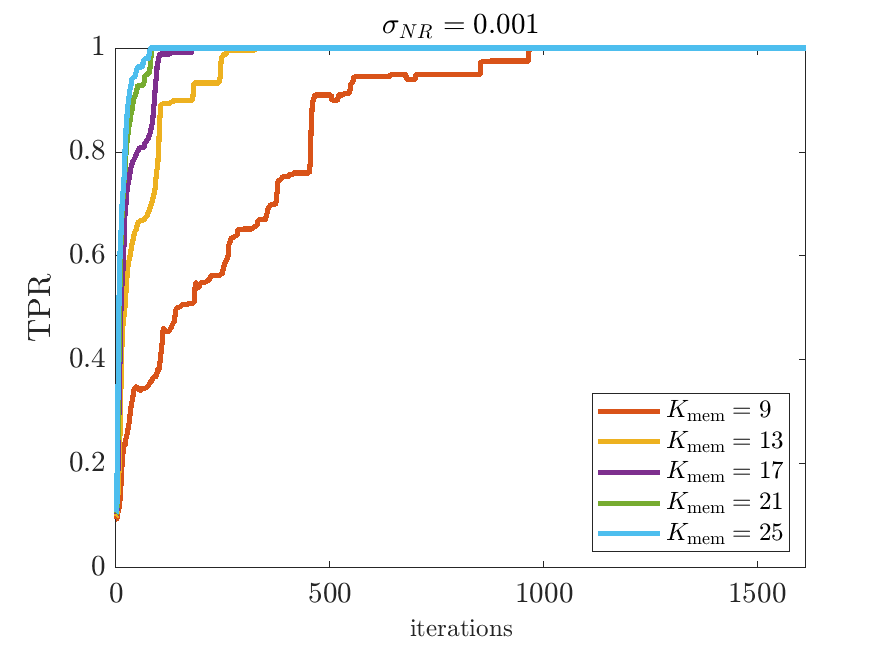} \\ 
\includegraphics[trim={0 0 20 0},clip,width=0.4\textwidth]{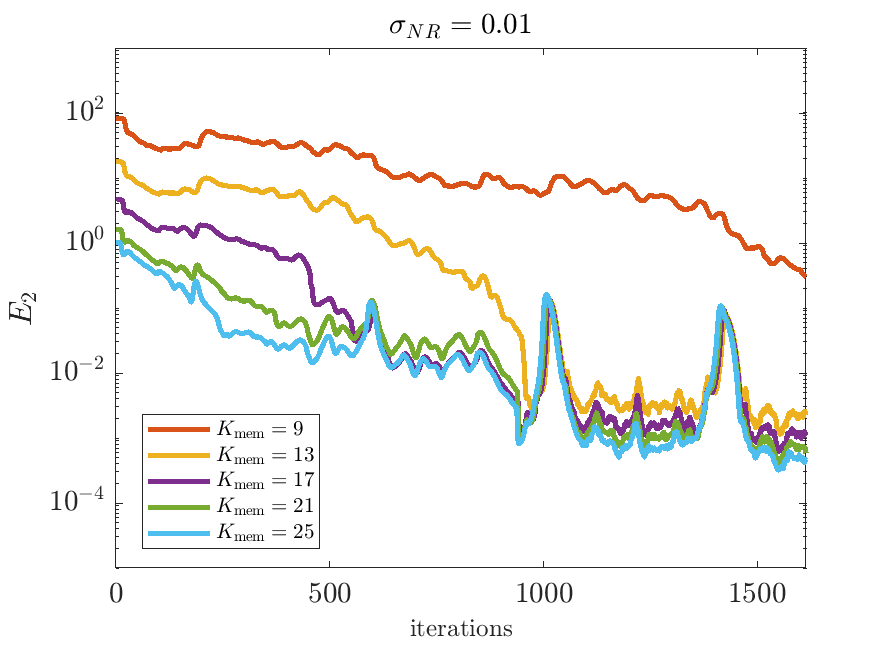} & 
\includegraphics[trim={0 0 20 0},clip,width=0.4\textwidth]{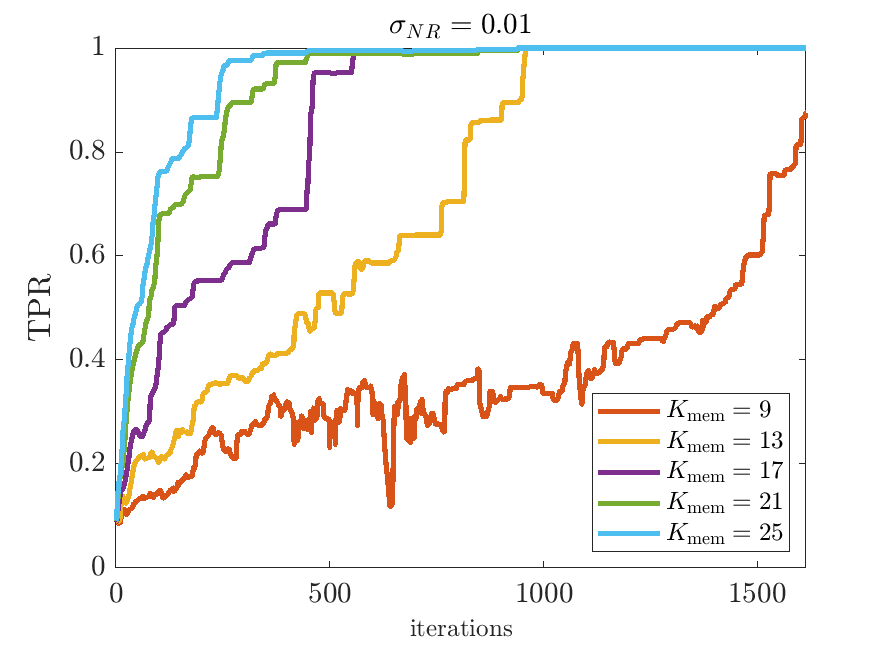} \\ 
\includegraphics[trim={0 0 20 0},clip,width=0.4\textwidth]{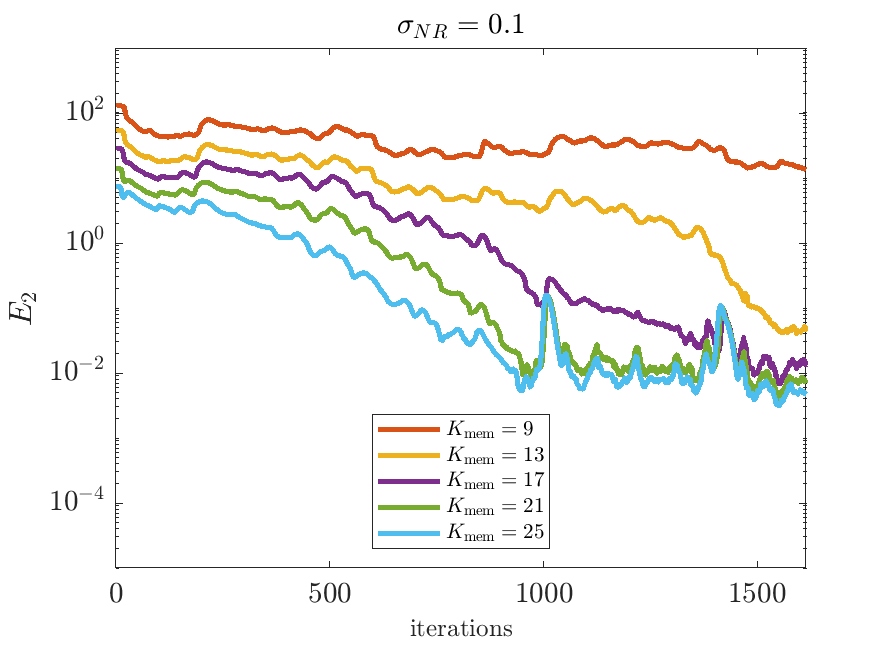} & 
\includegraphics[trim={0 0 20 0},clip,width=0.4\textwidth]{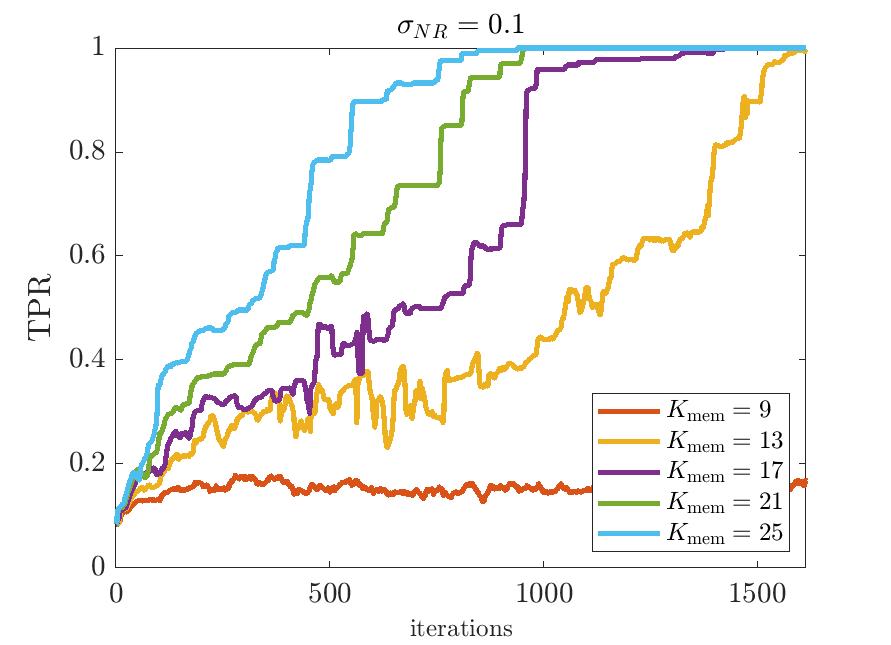} 
\end{tabular}
\caption{Online identification of the variable medium nonlinear wave equation \eqref{SGeq} for $K_\text{mem}\in\{9,13,17,21,25\}$ and (top to bottom) $\sigma_{NR}\in\{0,0.001,0.01,0.1\}$. Left: average coefficient error $E_2(\what^{(t)})$. Right: average total positivity ratio TPR$(\what^{(t)})$.}
\label{fig:W2D}
\end{figure}

Figure \ref{fig:W2D} shows robust recovery for $K_\text{mem}\in \{13,17,21,25\}$ up to $\sigma_{NR} =0.1$, with rapid identification for small noise. This is despite abrupt changes in the wavespeed $c$. For $K_\text{mem}=9$ we see recovery up to $\sigma_{NR}=0.001$, indicating that for larger noise 9 points in time is insufficient to discretize the integrals $\partial_{tt}\psi * u$ accurately, analogous to the case $K_\text{mem}=5$ for KS.

The left panel of Figure \ref{fig:W2D} shows that once the system is identified, abrupt changes in the wavespeed temporarily increase the coefficient error $E_2$, but the correct support $S^\star$ remains identified and the errors swiftly decay. In Figure \ref{fig:W2Dc} we plot the average learned wavespeed $\widehat{c}(t)$ as well as the maximum and minimum values of $\widehat{c}(t)$ attained over all 100 trials, revealing that increasing $K_\text{mem}$ from 17 to 25 leads to a significant decrease in the variance of $\widehat{c}$ after the system has been identified. This is purely an affect of using the weak form to discretize the time derivatives, and demonstrates that even under large noise and abruptly changing coefficients, the algorithm is able to maintain support recovery and accuracy.


\subsection{Wave equation in 3D}

\begin{equation}\label{W3D}
\partial_{tt} u = \partial_{xx}u+\partial_{yy}u+\partial_{zz}u
\end{equation}

For our last example we treat the linear wave equation in 3D. Exact data $\Ubf_{exact}$ has dimensions $128\times 128 \times 128\times  960$ in $(x,y,t)$ with resolution $(\Delta x,\Delta t)= (0.0491,0.0122)$. Each snapshot $\Ubf^{(t)}$ is $16.8$ Mb and online iterations take approximately $1.3$ seconds.

\begin{figure}
\centering
\begin{tabular}{cc}
\includegraphics[trim={0 0 20 0},clip,width=0.4\textwidth]{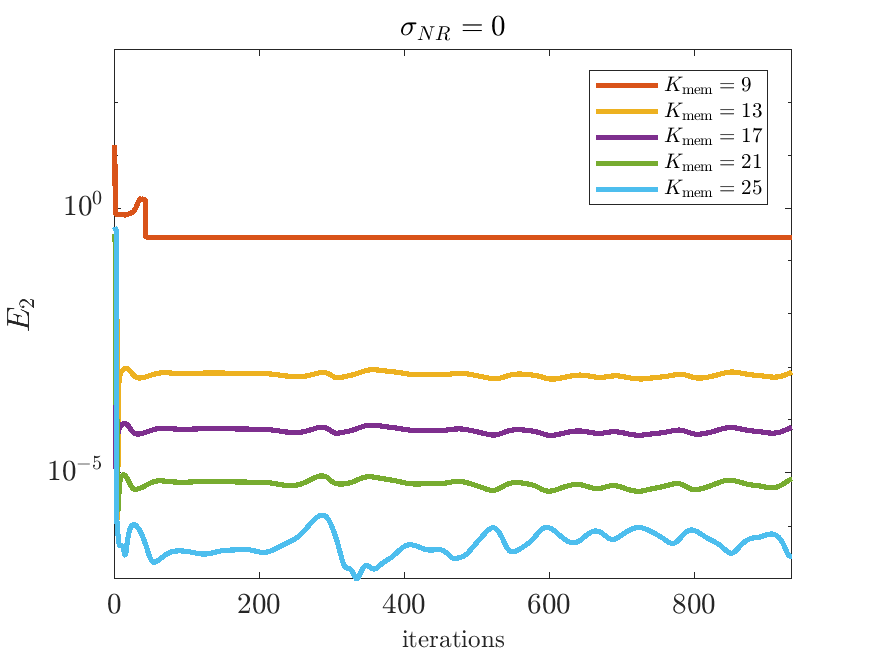} & 
\includegraphics[trim={0 0 20 0},clip,width=0.4\textwidth]{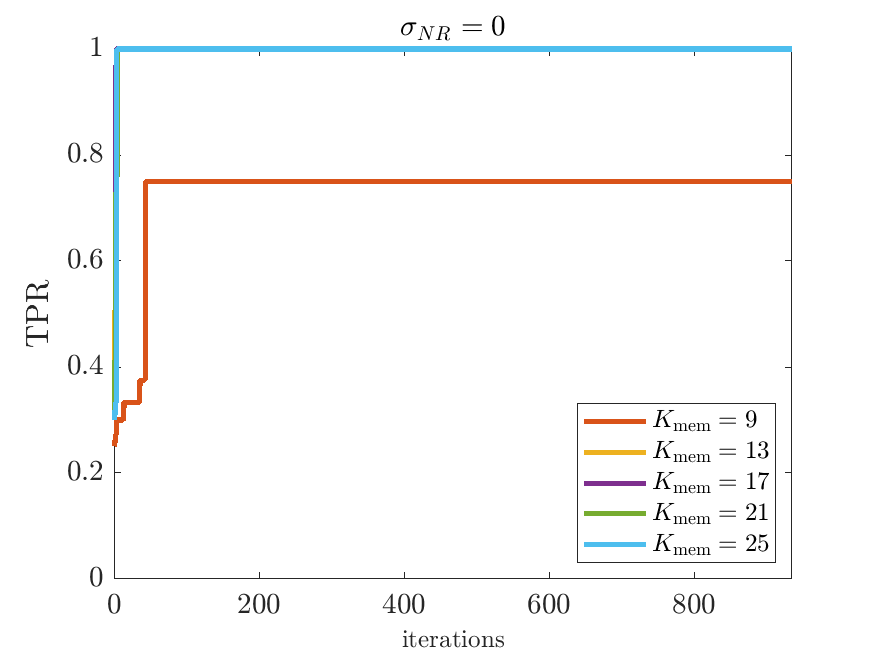} \\ 
\includegraphics[trim={0 0 20 0},clip,width=0.4\textwidth]{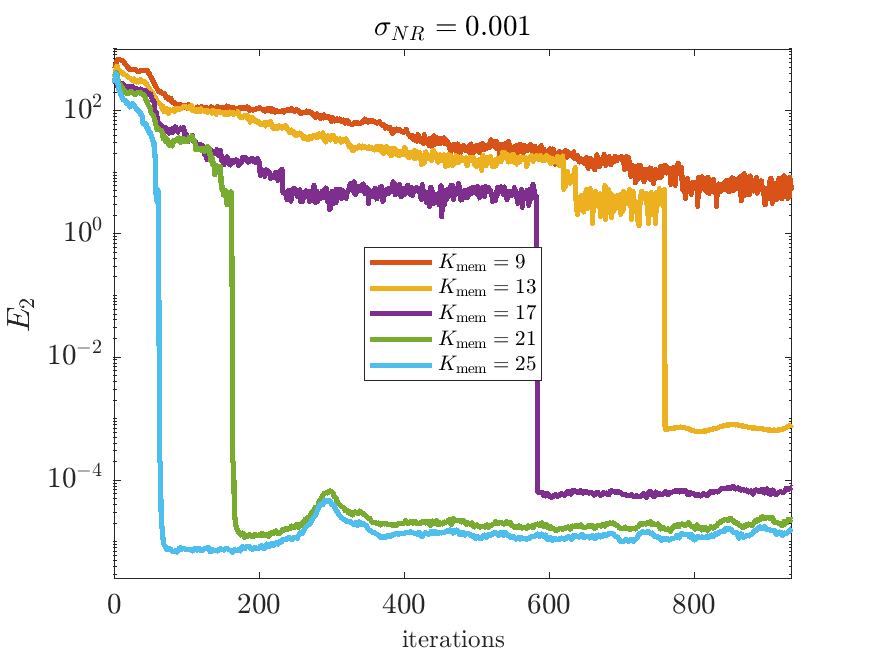} & 
\includegraphics[trim={0 0 20 0},clip,width=0.4\textwidth]{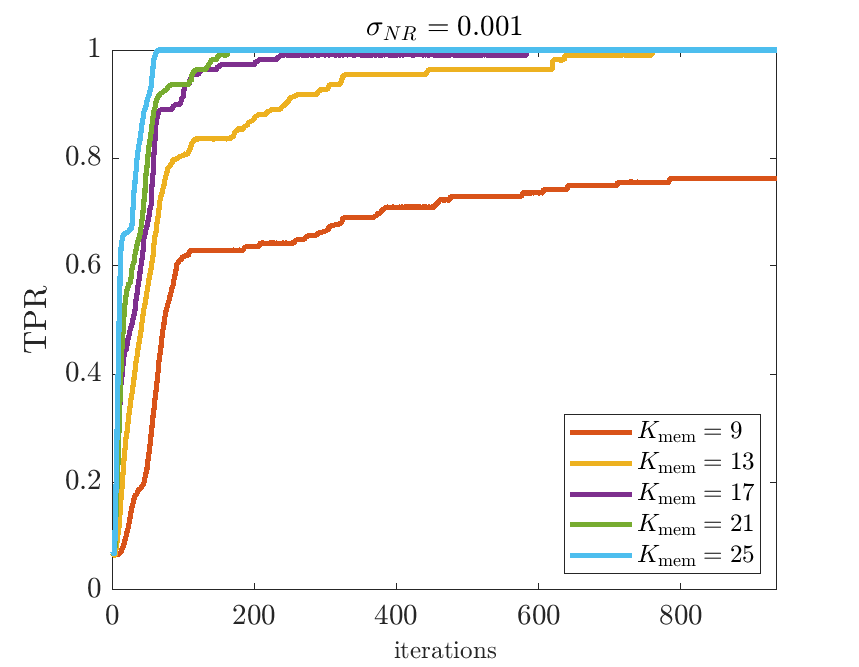} \\ 
\includegraphics[trim={0 0 20 0},clip,width=0.4\textwidth]{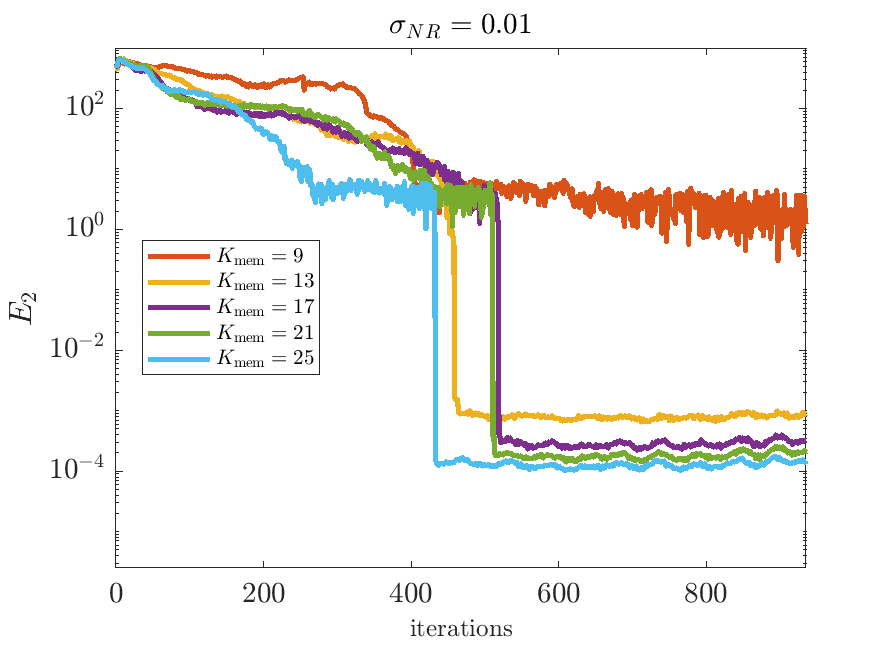} & 
\includegraphics[trim={0 0 20 0},clip,width=0.4\textwidth]{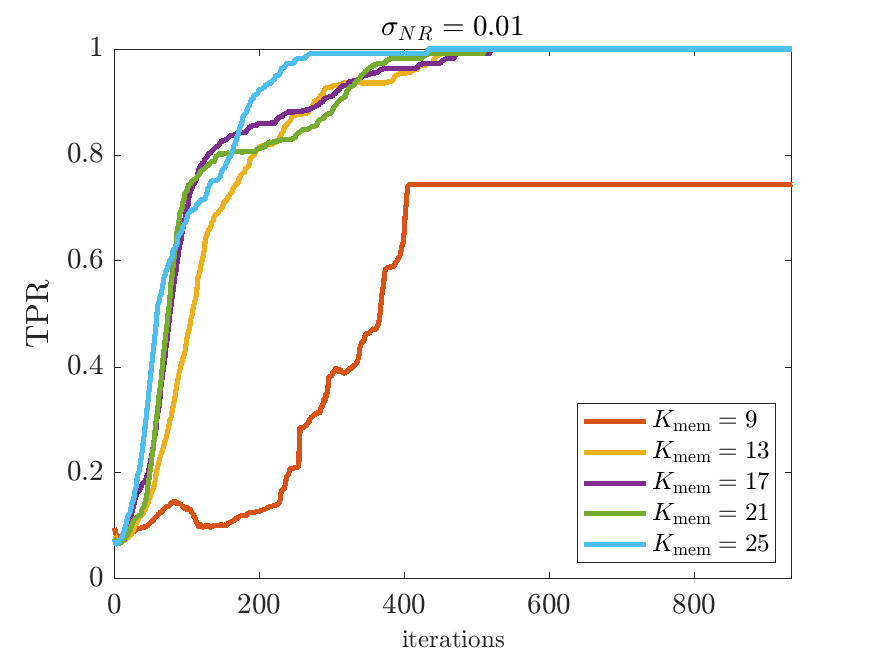} \\ 
\includegraphics[trim={0 0 20 0},clip,width=0.4\textwidth]{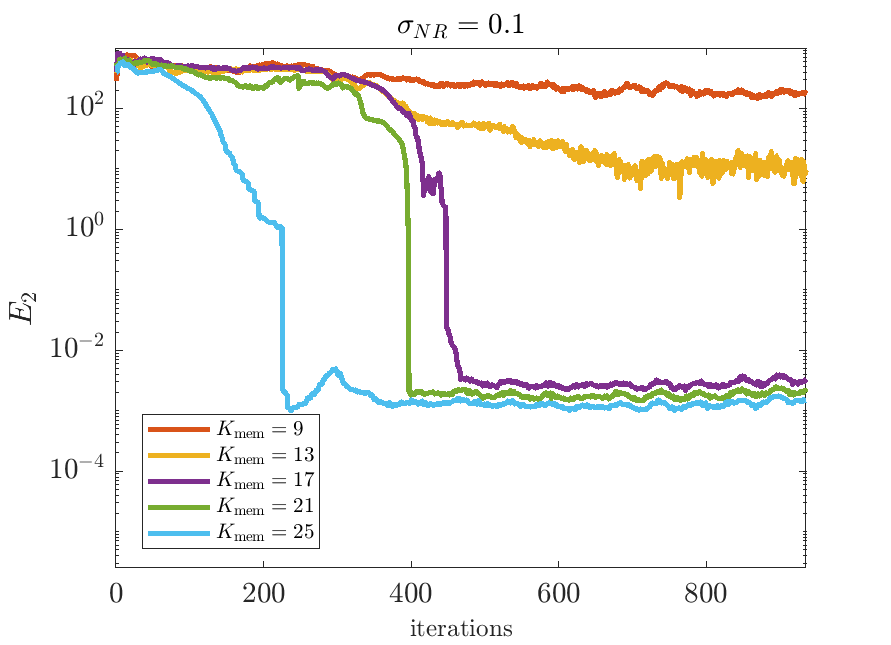} & 
\includegraphics[trim={0 0 20 0},clip,width=0.4\textwidth]{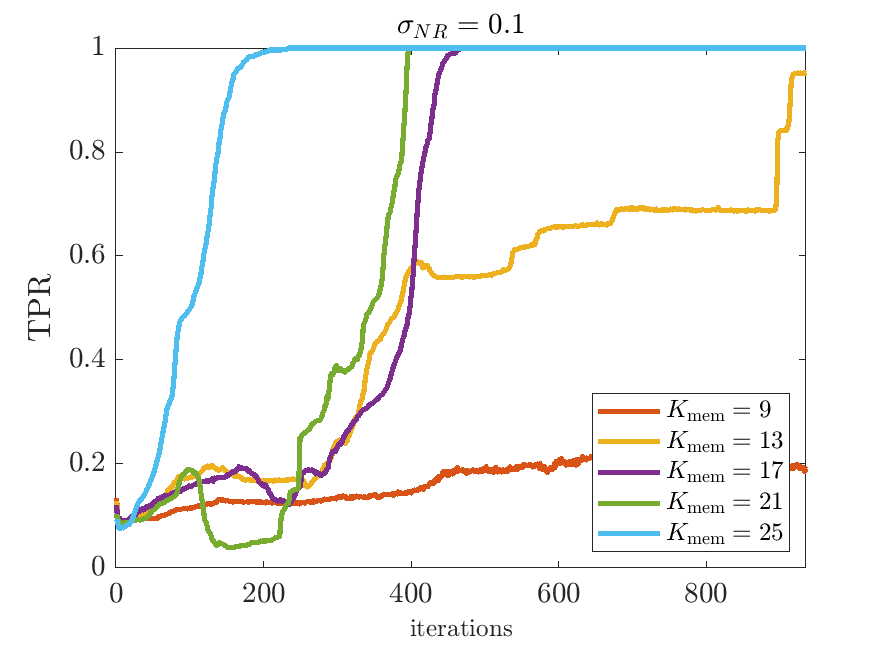}
\end{tabular}
\caption{Online identification of the wave equation in three spatial dimensions \eqref{W3D} for $K_\text{mem}\in\{9,13,17,21,25\}$ and (top to bottom) $\sigma_{NR}\in\{0,0.001,0.01,0.1\}$. Left: average coefficient error $E_2(\what^{(t)})$. Right: average total positivity ratio TPR$(\what^{(t)})$.}
\label{fig:W3D}
\end{figure}

Results are depicted in Figure \ref{fig:W3D}. We again find robust recovery for $K_\text{mem}\in \{13,17,21,25\}$ up to $\sigma_{NR}=0.1$, although in $5\%$ of trials at $\sigma_{NR}=0.1$ the $K_\text{mem}=13$ case finds a spurious term $\approx -0.8 u$. Even at $\sigma_{NR}=0.1$ the coefficients are accurate to more than 2 digits once recovered for $K_\text{mem}\geq 17$. For $K_\text{mem}=9$ we see poor performance for the same reason as above with W2D, but now manifesting as recovery of the spurious term $\approx -0.8u$, indicating that the inaccurate computation of $\partial_{tt}\psi * u$ produces spurious damping. This is not an altogether unreasonable affect if computing $\partial_{tt}\psi * u$ numerically is viewed as an attenuated second derivative calculation, although it does imply that higher-order time derivatives require more snapshots to be saved in memory.

\section{Conclusions}

We have demonstrated on several protoypical examples, and over a wide range of noise and memory scenarios, the viability of an online algorithm for PDE identification based on the weak-form sparse identification of nonlinear dynamics algorithm (WSINDy). The core of the method combines a weak-form discretization of candidate PDEs with the online proximal gradient descent algorithm applied directly to the least squares cost function with $\ell_0$-pseudo-norm regularization \eqref{WSINDyOL}. Compared with the more common approach of regularizing the $\ell_0$-pseudo-norm (e.g.\ with $\nrm{\cdot}_1$ or weighted variants \cite{candes2008enhancing}), we find that directly applying prox$_{\lambda \nrm{\cdot}_0}$, leading to hard thresholding, and adaptively selecting $\lambda_t$, exhibits good performance in efficiently identifying systems, handling noise, and tracking time-varying coefficients. 

Numerical experiments with an abruptly changing wavespeed indicate that our method is a lightweight counterpart to existing methods for variable coefficients (e.g. \cite{rudy2019data}), which may be of independent interest in the control of wave equations in variable-media (\cite{fante1971transmission,felsen1970wave,ning2010stabilization,seymour1987exact,chen1979control,vila2017bloch}). Examination of the wave equation in 3D also offers a different perspective on PDE identification in higher dimensions: problems with large datasets can be implemented in an online data-streaming fashion (not necessarily along the time axis as implemented here). It may therefore be advantageous from the standpoint of memory usage to solve certain batch problems in the online manner we have presented.

The algorithm's successes warrant further investigation in a number of areas. While we have characterized stationary points of the batch algorithm and proved boundedness of the average dynamic regret, we leave a more complete analysis to future work. In particular, one could analyze the error $\nrm{\wbf^{(t)}-\wstar(t)}$ as a function of library $\Theta$, test function $\psi$, data sampling rates $(\Delta x,\Delta t)$, memory size $K_\text{mem}$, noise ratio $\sigma_{NR}$, etc. It may also advantageous to design adaptive schemes which update $\Theta$ and $\psi$ throughout the course of the algorithm, depending on the dynamics of the data and previously learned equations. Nevertheless, the current framework is well-suited for a large variety of problems and opens the door to online PDE identification as well as the possibility of solving batch problems in an online manner.

\section{Acknowledgements}

his research was supported in part by the NSF/NIH Joint DMS/NIGMS Mathematical Biology Initiative grant R01GM126559, in part by the NSF Mathematical Biology MODULUS grant 2054085, and in part by the NSF Computing and Communications Foundations grant 2054085. This work also utilized resources from the University of Colorado Boulder Research Computing Group, which is supported by the National Science Foundation (awards ACI-1532235 and ACI-1532236), the University of Colorado Boulder, and Colorado State University.

\bibliographystyle{plain}
\bibliography{researchCU}

\appendix

\section{Column scaling and non-uniform thresholds}\label{App:tech}

For stability, we normalize the columns of $\Gbf^{(t)}$ at each step, defining $\tilde{\Gbf}^{(t)} = \Gbf^{(t)}
\Mbf^{(t)}$ with 
\[\Mbf^{(t)} = \text{diag}\left(\nrm{\Gbf^{(t)}_{1}}_2^{-1},\dots,\nrm{\Gbf^{(t)}_{IJ}}_2^{-1}\right).\] 
In particular, this allows for a larger stepsize $\tilde{\alpha}_t =1/\nrm{(\tilde{\Gbf}^{(t)})^T\tilde{\Gbf}^{(t)}_{S_t}}_2$ and leads to a reasonable estimate $\tilde{\alpha}_t = 1/\sqrt{|S_t|IJ}$ for a stepsize that does not require computation of the matrix $2$-norm.

For more flexibility, we allow for non-uniform thresholding. For a set of thresholds $\pmb{\lambda}\in \Rbb^{IJ}$, we define the non-uniform thresholding operator $H_{\pmb{\lambda}}$ by
\[(H_{\pmb{\lambda}}(\xbf))_i = \begin{dcases} \xbf_i, & |\xbf_i|\geq \pmb{\lambda}_i \\ 
0, & \text{otherwise}.\end{dcases}\] 
This happens to be the proximal operator of the non-uniform $\ell_0$-norm
\begin{equation}
\label{app:ell0}
\nrm{\xbf}_{0,\pmb{\lambda}}:=\sum_{i=1}^{IJ} \pmb{\lambda}_i^2 \textbf{1}_{\Rbb\setminus\{0\}}(\xbf_i),
\end{equation}
where $\nrm{\xbf}_{0,\pmb{\lambda}} = \lambda^2\nrm{\xbf}_0$ when $\pmb{\lambda} = (\lambda,\dots,\lambda)$. The resulting online cost function being minimized after incorporation of both non-uniform thresholding and column rescaling is 
\begin{equation}\label{app:F1}
\tilde{F}_t(\wbf;\pmb{\lambda}^{(t)}) = \frac{1}{2}\nrm{\tilde{\Gbf}^{(t)}\wbf-\bbf^{(t)}}_2^2+\frac{1}{2}\nrm{\wbf}_{0,(\Mbf^{(t)})^{-1}\pmb{\lambda}^{(t)}},
\end{equation}
whose fixed points $\tilde{\wbf}^{\star,(t)}$ coincide with those of the desired cost function
\begin{equation}\label{app:F2}
F_t(\wbf;\pmb{\lambda}^{(t)}) = \frac{1}{2}\nrm{\Gbf^{(t)}\wbf-\bbf^{(t)}}_2^2+\frac{1}{2}\nrm{\wbf}_{0,\pmb{\lambda}^{(t)}}
\end{equation}
after a diagonal transformation $\wbf^{\star,(t)} = \Mbf^{(t)}\tilde{\wbf}^{\star,(t)}$. With these two pieces, the online algorithm for \eqref{app:F1} becomes 
\[\begin{dcases}
\tilde{\wbf}^{(t)} = (\Mbf^{(t)})^{-1}\widehat{\wbf}^{(t)} \\ 
\what_{t+\Delta t} = H_{\tilde{\alpha}_t\pmb{\lambda}^{(t)}}\left(\Mbf^{(t)}\left(\tilde{\wbf}^{(t)}-\tilde{\alpha}_t(\tilde{\Gbf}^{(t)})^T\left(\tilde{\Gbf}^{(t)}\tilde{\wbf}^{(t)}-\bbf^{(t)}\right)\right)\right),
\end{dcases}\]
however this can equivalently be written in terms of the desired coefficients $\what^{(t)}$ as
\begin{equation}
\what^{(t+1)} = H_{\tilde{\alpha}_t\pmb{\lambda}^{(t)}}\left(\what^{(t)}-\tilde{\alpha}_t(\Mbf^{(t)})^2(\Gbf^{(t)})^T\left(\Gbf^{(t)}\what^{(t)}-\bbf^{(t)}\right)\right).
\end{equation}
In direct analogy to the batch WSINDy thresholding scheme \eqref{MSTLS1}-\eqref{MSTLSbnds}, we use thresholds \[\pmb{\lambda}^{(t)}=\max\left(1,\nrm{\bbf^{(t)}}\text{diag}(\Mbf^{(t)})\right)\lambda_t,\] 
which eliminate small coefficient values $\min_{i\in S_t}|\widehat{\wbf}^{(t)}_i|\geq \lambda_t$ as well as small terms in the sense of dominant balance with respect to $\bbf^{(t)}$: 
\[\min_{i\in S_t}\frac{\nrm{\Gbf^{(t)}_{i}\widehat{\wbf}^{(t)}_i}_2}{\nrm{\bbf^{(t)}}_2} \geq \lambda_t.\]
The update rule \eqref{lambdaupdate} for $\lambda_t$ is unchanged after replacing $F_t(\wbf;\lambda_t)$ with $F_t(\wbf;\pmb{\lambda}^{(t)})$ defined in \eqref{app:F2}.
\section{Implementation and Computational Complexity}

The offline phase has four components:
\begin{enumerate}
\item Initialize hyperparameters $\psi(\xbf,t)=\phi(\xbf)\theta(t)$, $\Theta=\{D^{\pmb{\alpha}^{(i)}}f_j\}_{i=0,j=1}^{I,J}$, $\Delta \lambda$, $\lambda_{\max}$, $\lambda_0$, where the test function $\psi$ is either prescribed manually or selected using the changepoint algorithm from \cite{MessengerBortz2021JComputPhys} using the initial $K_\text{mem}$ slices $\{\Ubf^{(0)},\dots,\Ubf^{((K_\text{mem}-1)\Delta t)}\}$.
\item Compute and store the Fourier transforms $\{\widehat{D^{\pmb{\alpha}^{(i)}}\psi}\}_{i=0}^I$ to reuse at each step when computing convolutions (recall $\psi$ is separable so this storage cost is negligible). 
\item Compute initial library of spatially integrated terms 
\[\Psi\ :=\ \{\Psi^{(t)}\}_{t=0}^{(K_\text{mem}-1)\Delta t}\ := \ \left\{\left\{D^{\pmb{\beta}^{(i)}}\phi*f_j(\Ubf^{(t)})(\CalQ_\xbf,t)\right\}_{i=0,j=1}^{I,J}\right\}_{t=0}^{(K_\text{mem}-1)\Delta t}\]
where $\pmb{\beta}^{(i)} = (\pmb{\alpha}^{(i)}_1,\dots,\pmb{\alpha}^{(i)}_d)$ is the spatial part of the multi-index $\pmb{\alpha}^{(i)}$ operating on the spatial part  $\phi$ of the test function $\psi$ (recall that $\CalQ_\xbf$ is the set of spatial points over which convolutions are evaluated, also equal to the number of rows in $\Gbf^{(t)}$). 
\item Compute initial weights $\what^{(0)} = (\Gbf^{(0)})^\dagger \bbf^{(0)}$ where $\bbf^{(0)}$ and $\Gbf^{(0)}$ are obtained by integrating the elements of $\Psi$ in time against the corresponding temporal test functions $D^{\pmb{\alpha}^{(i)}_{d+1}}\theta$.
\end{enumerate}
For each $t$ in the online phase we compute $\Psi^{(t)}$ using only the incoming slice $\Ubf^{(t)}$, which replaces $\Psi^{(t-K_\text{mem}\Delta t)}$ in memory. $(\Gbf^{(t)},\bbf^{(t)})$ are then computed by integrating the elements of $\Psi$ in time against the corresponding temporal test functions $D^{\pmb{\alpha}^{(i)}_{d+1}}\theta$, which amounts to a series of dot products between length-$K_\text{mem}$ vectors. Computation of $\Gbf^{(t)}$ at each time $t$ thus requires $J|\Xbf|$ function evaluations $f_j(\Ubf^{(t)})$  (each counted as 1 floating point operation (flop)) followed by $IJ$ convolutions against $D^{\pmb{\beta}^{(i)}}\phi$, and finally integration in time. The total flop count at each step is at most 
\[J|\Xbf|\left(1\ +\ CI\log N\ +\ 2IK_{\text{mem}} \frac{|\CalQ_\Xbf|}{|\Xbf|}\right)\]
where $C$ is such that $\xbf*\ybf$ costs $CN \log N$ using FFTs for length-$N$ vectors $\xbf$ and $\ybf$, minus the cost of one FFT (since we have precomputed these for $D^{\pmb{\beta}^{(i)}}\phi$) and $N \approx |\Xbf|^{1/d}$ is the one-dimensional length scale of the data. In other words, only 
\[F = J\left(1+CI\log N +IK_\text{mem}\frac{|\CalQ_\xbf|}{|\Xbf|}\right)\]
flops are performed per incoming data point in $\Ubf^{(t)}$ (and a more careful analysis leads to a lower cost in the factor $CI\log N$ by incorporating the subsampling $\Xbf\to \CalQ_\xbf$). Note that $F$ does not depend on the spatial dimension $d$ of the data set (except through library term $I$, which might increase with $d$ as more differential operators become added). The total working memory $W$ to store $\Psi$ and ($\Gbf^{(t)}$, $\bbf^{(t)}$) as outlined above is given by 
$W = IJ|\CalQ_\xbf|K_\text{mem}+ (I+1)J|\CalQ_\xbf|$ double-precision floating point numbers (DPs).

\begin{rmrk}
 There are several natural choices to consider to either decrease storage restrictions or increase computational speed. However, it is not clear that the anticipated savings will manifest. For instance, we could instead store the spatial Fourier transforms of the nonlinearities $\{\widehat{f_j(\Ubf^{(t)})}\}_{j=1,t=\ell\Delta t}^{J, (\ell+K_\text{mem}-1)\Delta t}$, resulting in a working memory of $J\cdot K_\text{mem}\cdot |\Xbf|$ instead of $IJ|\CalQ_\Xbf|K_\text{mem}$ to store $\Psi$. This would require that we compute spatial convolutions over all $K_\text{mem}$ time slices at each time point, instead of spatial convolutions over just the incoming time slice $\Ubf^{(t)}$, hence resulting in a $K_\text{mem}$-fold increase in computation time, as this is the leading-order cost. In addition, the storage ``savings'' may actually be worse, specifically if $I|\CalQ_\xbf|\leq |\Xbf|$. We believe that the method outlined above provides a near-optimal balance of computational complexity and storage requirements, with a heavier emphasis on reducing computational complexity. 

\section{Proof of Lemma \ref{lemm}}\label{app:lemm}
 
Consider $\wbf$ such that one of the following holds:
\begin{enumerate}[label=(\roman*)]
\item $\wbf$ is a local minimizer of \eqref{sparserec}
\item $\wbf = H_\lambda\left(\wbf - \Gbf^T\left(\Gbf \wbf - \bbf\right)\right)$
\item With $S=\supp{\wbf}$, we have that $\wbf_S \in \argmin_\zbf\nrm{\Gbf_{S}\zbf-\bbf}_2^2$ and
\[\max_{i\in S^c}\left\vert \Gbf^T_{i}(\Gbf\wbf-\bbf) \right\vert<\lambda \leq \min_{i\in S} \left\vert \wbf_i\right\vert.\]
\end{enumerate}
Then it holds that $(ii)\iff(iii) \implies (i)$. Moreover, if $\wbf$ a global minimizer, then $(i)\implies (iii)$.\\

\begin{proof}
$(iii)\implies (ii)$ is immediate. To show $(ii)\implies (iii)$, let $S=\supp{\wbf}$. Then we have  
\[\wbf_S = \wbf_S-\Gbf^T_{S}(\Gbf\wbf-\bbf),\]
which implies that $\min_{i\in S}|\wbf_i|\geq \lambda$ so that $ \Gbf^T_{S}\Gbf_{S}\wbf_S = \Gbf^T_{S}\bbf$, so that $\wbf_S \in \argmin_\zbf\nrm{\Gbf_{S}\zbf-\bbf}_2^2$. On $S^c$ we have
\[ H_\lambda\left(\Gbf_{S^c}^T\left(\Gbf\wbf -\bbf\right)\right)=0\implies \max_{i\in S^c}\left\vert \Gbf^T_{i}(\Gbf\wbf-\bbf)\right\vert<\lambda.\]
To show that $(ii)$ and $(iii)$ imply $(i)$, we note that under usual assumptions of two closed, convex and proper functions $f$ and $g$, we have
\[\wbf\in \text{prox}_g\left(\wbf - \partial f(\wbf)\right) \iff 0\in \partial f(\wbf)+\partial g(\wbf) \implies \wbf\in \text{argmin}(f+g),\]
however $\nrm{\cdot}_0$ is clearly not convex\footnote{In fact the subdifferential $\partial\nrm{\cdot}_0(\wbf)=\emptyset$ unless $\wbf=\textbf{0}$, upon which $\partial\nrm{\cdot}_0(\wbf)=\{\textbf{0}\}$.}. Instead we can directly show that for a perturbed vector $\tilde{\wbf} = \wbf +\pmb{\eta}$, for suitably small $\nrm{\pmb{\eta}}$ the objective is non-decreasing. Using that $\wbf_S \in \argmin_\zbf\nrm{\Gbf_{S}\zbf-\bbf}_2^2$, let $\Pbf_S^\perp$ be the projection onto $\{\text{span}(\Gbf_{S})\}^\perp$. The difference in objective $F$ is then given by 
\[F(\tilde{\wbf};\lambda)-F(\wbf;\lambda) = \frac{1}{2}\left(\nrm{\Pbf_S^\perp\bbf+\Gbf\pmb{\eta}}_2^2-\nrm{\Pbf_S^\perp\bbf}_2^2\right)+\frac{\lambda^2}{2}\left(\nrm{\tilde{\wbf}}_0-\nrm{\wbf}_0\right)\]
\[=\frac{1}{2}\nrm{\Gbf\pmb{\eta}}_2^2+\lan \Pbf_S^\perp \bbf,\Gbf\pmb{\eta}\ran +  \frac{\lambda^2}{2}\left(\nrm{\tilde{\wbf}}_0-\nrm{\wbf}_0\right).\]
If $\supp{\pmb{\eta}}\subset \supp{\wbf}$ and $\nrm{\pmb{\eta}}_\infty<\lambda$, then $\nrm{\tilde{\wbf}}_0=\nrm{\wbf}_0$ and $\lan \Pbf_S^\perp \bbf,\Gbf\pmb{\eta}\ran=0$, hence $F(\tilde{\wbf};\lambda)-F(\wbf;\lambda) \geq 0$, with equality only if $\Gbf\pmb{\eta}=\mathbf{0}$, which in particular is not possible when $\Gbf_{S}$ is full rank unless $\pmb{\eta}=0$. If $\supp{\pmb{\eta}}\notin S$ and $\nrm{\pmb{\eta}}_\infty<\lambda$, then $\Pbf^\perp_S\bbf=\mathbf{0}$ implies a strict increase in $F$, while if $\Pbf^\perp_S \bbf\neq \mathbf{0}$ then
\[\nrm{\pmb{\eta}_{S^c}}_2 < \ep:=\frac{\lambda^2}{2}\frac{1}{\nrm{\Pbf^\perp_S\bbf}_2\nrm{\Gbf_{S^c}}_2},\]
implies a strict increase in $F$. To see this, note that $\lan \Pbf_S^\perp,\Gbf\pmb{\eta}\ran = \lan \Pbf_S^\perp, \Gbf_{S^c}\pmb{\eta}_{S^c}\ran$ implies the bound
\[F(\tilde{\wbf};\lambda)-F(\wbf;\lambda) \geq -\nrm{\Pbf^\perp_S\bbf}_2\nrm{\Gbf_{S^c}}_2\nrm{\pmb{\eta}_{S^c}}_2+\frac{\lambda^2}{2}>0.\]
Note that $\ep$ is not tight. Combining these conditions gives a ball around $\wbf$ over which $F$ is non-decreasing, hence $\wbf$ is a local min. Finally, that $\wbf$ a global minimizer implies $(iii)$ can be found in \cite{zhang2019convergence}.
\end{proof}
 
\end{rmrk}
 



\end{document}

%% file: ex_shared.tex
\usepackage{lipsum}
\usepackage{amsfonts}
\usepackage{graphicx}
\usepackage{epstopdf}
\usepackage{subcaption}
\usepackage{algorithmic}
\ifpdf
  \DeclareGraphicsExtensions{.eps,.pdf,.png,.jpg}
\else
  \DeclareGraphicsExtensions{.eps}
\fi

\newsiamremark{remark}{Remark}
\newsiamremark{hypothesis}{Hypothesis}
\crefname{hypothesis}{Hypothesis}{Hypotheses}
\newsiamthm{claim}{Claim}

\headers{Online Weak-form Sparse Identification}{D. Messenger, E. Dall'Anese, D. Bortz}

\title{Online Weak-form Sparse Identification \\ of Partial Differential Equations}

\author{Daniel Messenger\thanks{Department of Applied Mathematics, University of Colorado, Boulder, CO 80309-0526, USA. 
	(\email{daniel.messenger@colorado.edu}, \email{david.bortz@colorado.edu}).}
\and Emiliano Dall'Anese\thanks{Department of Department of Electrical, Computer, and Energy Engineering, University of Colorado, Boulder, CO 80309-0526, USA. 
	(\email{emiliano.dallanese@colorado.edu}).}
\and David Bortz\footnotemark[1]}

\usepackage{amsopn}

\usepackage[latin9]{inputenc}
\usepackage{geometry}
\usepackage{xcolor}
\usepackage{pdfcolmk}
\usepackage{float}
\usepackage{mathtools}
\usepackage{amsbsy}
\usepackage{amstext}
\usepackage{amssymb}
\usepackage{graphicx}
\usepackage{
enumitem}
\PassOptionsToPackage{normalem}{ulem}
\usepackage{ulem}

\graphicspath{{figures/}}

\makeatletter

\providecolor{lyxadded}{rgb}{0,0,1}
\providecolor{lyxdeleted}{rgb}{1,0,0}

\DeclareRobustCommand{\lyxsout}[1]{\ifx\\#1\else\sout{#1}\fi}

\numberwithin{equation}{section}
\numberwithin{figure}{section}
\theoremstyle{plain}
\newtheorem{thm}{\protect\theoremname}
  \theoremstyle{remark}
  
  \theoremstyle{plain}

\input{notation}

\usepackage{stmaryrd}

\makeatother

  \providecommand{\lemmaname}{Lemma}
  \providecommand{\remarkname}{Remark}
\providecommand{\theoremname}{Theorem}


%% file: notation.tex
\newcommand{\Cbb}{\mathbb{C}}

\newcommand{\Nbb}{\mathbb{N}}

\newcommand{\Rbb}{\mathbb{R}}

\newcommand{\bbf}{\mathbf{b}}

\newcommand{\Gbf}{\mathbf{G}}

\newcommand{\Ibf}{\mathbf{I}}

\newcommand{\Mbf}{\mathbf{M}}

\newcommand{\Pbf}{\mathbf{P}}

\newcommand{\Ubf}{\mathbf{U}}

\newcommand{\wbf}{\mathbf{w}}

\newcommand{\xbf}{\mathbf{x}}
\newcommand{\Xbf}{\mathbf{X}}
\newcommand{\ybf}{\mathbf{y}}

\newcommand{\zbf}{\mathbf{z}}

\newcommand{\ep}{\epsilon}

\newcommand{\CalI}{{\mathcal{I}}}

\newcommand{\CalL}{{\mathcal{L}}}

\newcommand{\CalQ}{{\mathcal{Q}}}

\newcommand{\CalT}{{\mathcal{T}}}

\newcommand{\nrm}[1]{\left\Vert {#1} \right\Vert}

\newcommand{\supp}[1]{\text{supp}\left(#1\right)}

\renewcommand{\tilde}[1]{\widetilde{#1}}

\newcommand{\argmin}{\text{argmin}}

\newcommand{\lan}{\left\langle}
\newcommand{\ran}{\right\rangle}

\newcommand{\wstar}{\wbf^\star}
\newcommand{\what}{{\widehat{\wbf}}}

\newtheorem{lemm}{Lemma}[section]

\theoremstyle{definition}
\newtheorem{rmrk}{Remark}[section]


%% file: MSML2022_arXiv.bbl
\begin{thebibliography}{10}

\bibitem{bakarji2021data}
Joseph Bakarji and Daniel~M Tartakovsky.
\newblock Data-driven discovery of coarse-grained equations.
\newblock {\em Journal of Computational Physics}, 434:110219, 2021.

\bibitem{boninsegna2018sparse}
Lorenzo Boninsegna, Feliks N{\"u}ske, and Cecilia Clementi.
\newblock Sparse learning of stochastic dynamical equations.
\newblock {\em The Journal of chemical physics}, 148(24):241723, 2018.

\bibitem{brunton2016discovering}
Steven~L Brunton, Joshua~L Proctor, and J~Nathan Kutz.
\newblock Discovering governing equations from data by sparse identification of
  nonlinear dynamical systems.
\newblock {\em Proceedings of the national academy of sciences},
  113(15):3932--3937, 2016.

\bibitem{candes2008enhancing}
Emmanuel~J Candes, Michael~B Wakin, and Stephen~P Boyd.
\newblock Enhancing sparsity by reweighted $\ell_1$ minimization.
\newblock {\em Journal of Fourier analysis and applications}, 14(5):877--905,
  2008.

\bibitem{chen1979control}
Goong Chen.
\newblock Control and stabilization for the wave equation in a bounded domain.
\newblock {\em SIAM Journal on Control and Optimization}, 17(1):66--81, 1979.

\bibitem{dall2020optimization}
Emiliano Dall'Anese, Andrea Simonetto, Stephen Becker, and Liam Madden.
\newblock Optimization and learning with information streams: Time-varying
  algorithms and applications.
\newblock {\em IEEE Signal Processing Magazine}, 37(3):71--83, 2020.

\bibitem{dixit2019online}
Rishabh Dixit, Amrit~Singh Bedi, Ruchi Tripathi, and Ketan Rajawat.
\newblock Online learning with inexact proximal online gradient descent
  algorithms.
\newblock {\em IEEE Transactions on Signal Processing}, 67(5):1338--1352, 2019.

\bibitem{fante1971transmission}
R~Fante.
\newblock Transmission of electromagnetic waves into time-varying media.
\newblock {\em IEEE Transactions on Antennas and Propagation}, 19(3):417--424,
  1971.

\bibitem{fattahi2019learning}
Salar Fattahi, Nikolai Matni, and Somayeh Sojoudi.
\newblock Learning sparse dynamical systems from a single sample trajectory.
\newblock In {\em 2019 IEEE 58th Conference on Decision and Control (CDC)},
  pages 2682--2689. IEEE, 2019.

\bibitem{felsen1970wave}
L~Felsen and G~Whitman.
\newblock Wave propagation in time-varying media.
\newblock {\em IEEE Transactions on Antennas and Propagation}, 18(2):242--253,
  1970.

\bibitem{foster2020learning}
Dylan Foster, Tuhin Sarkar, and Alexander Rakhlin.
\newblock Learning nonlinear dynamical systems from a single trajectory.
\newblock In {\em Learning for Dynamics and Control}, pages 851--861. PMLR,
  2020.

\bibitem{10.5555/2526243}
Simon Foucart and Holger Rauhut.
\newblock {\em A Mathematical Introduction to Compressive Sensing}.
\newblock Birkh\"{a}user Basel, 2013.

\bibitem{hastie2009elements}
Trevor Hastie, Robert Tibshirani, Jerome~H Friedman, and Jerome~H Friedman.
\newblock {\em The elements of statistical learning: data mining, inference,
  and prediction}, volume~2.
\newblock Springer, 2009.

\bibitem{hazan2006efficient}
Elad Hazan.
\newblock {\em Efficient algorithms for online convex optimization and their
  applications}.
\newblock Princeton University, 2006.

\bibitem{hazan2009stochastic}
Elad Hazan and Satyen Kale.
\newblock On stochastic and worst-case models for investing.
\newblock {\em Advances in Neural Information Processing Systems}, 22, 2009.

\bibitem{hoi_online_2021}
Steven~C.H. Hoi, Doyen Sahoo, Jing Lu, and Peilin Zhao.
\newblock Online learning: {A} comprehensive survey.
\newblock {\em Neurocomputing}, 459:249--289, October 2021.

\bibitem{jialei_wang_online_2014}
{Jialei Wang}, {Peilin Zhao}, Steven C.~H. Hoi, and {Rong Jin}.
\newblock Online {Feature} {Selection} and {Its} {Applications}.
\newblock {\em IEEE Transactions on Knowledge and Data Engineering},
  26(3):698--710, March 2014.

\bibitem{kaiser2018sparse}
Eurika Kaiser, J~Nathan Kutz, and Steven~L Brunton.
\newblock Sparse identification of nonlinear dynamics for model predictive
  control in the low-data limit.
\newblock {\em Proceedings of the Royal Society A}, 474(2219):20180335, 2018.

\bibitem{koller2018learning}
Torsten Koller, Felix Berkenkamp, Matteo Turchetta, and Andreas Krause.
\newblock Learning-based model predictive control for safe exploration.
\newblock In {\em 2018 IEEE conference on decision and control (CDC)}, pages
  6059--6066. IEEE, 2018.

\bibitem{kopsinis_online_2011}
Yannis Kopsinis, Konstantinos Slavakis, and Sergios Theodoridis.
\newblock Online {Sparse} {System} {Identification} and {Signal}
  {Reconstruction} {Using} {Projections} {Onto} {Weighted} $\ell_1$ {Balls}.
\newblock {\em IEEE Transactions on Signal Processing}, 59(3):936--952, March
  2011.

\bibitem{ljung1999system}
L.~Ljung.
\newblock {\em System identification: theory for the user}.
\newblock 2nd edition Prentice-Hall, Upper Saddle River, NJ, 1999.

\bibitem{messenger2021learning}
Daniel~A Messenger and David~M Bortz.
\newblock Learning mean-field equations from particle data using {WSINDy}.
\newblock {\em arXiv preprint arXiv:2110.07756 (in revision at Physica D:
  Nonlinear Phenomena)}, 2021.

\bibitem{MessengerBortz2021JComputPhys}
Daniel~A. Messenger and David~M. Bortz.
\newblock Weak {{SINDy For Partial Differential Equations}}.
\newblock {\em J. Comput. Phys.}, 443:110525, October 2021.

\bibitem{MessengerBortz2021SIAMMultiscaleModelSimul}
Daniel~A. Messenger and David~M. Bortz.
\newblock Weak {{SINDy}}: {{Galerkin-Based Data-Driven Model Selection}}.
\newblock {\em SIAM Multiscale Model. Simul.}, 19(3):1474--1497, 2021.

\bibitem{nikolova2013description}
Mila Nikolova.
\newblock Description of the minimizers of least squares regularized with
  $\ell_0$-norm. uniqueness of the global minimizer.
\newblock {\em SIAM Journal on Imaging Sciences}, 6(2):904--937, 2013.

\bibitem{ning2010stabilization}
Zhen-Hu Ning and Qing-Xu Yan.
\newblock Stabilization of the wave equation with variable coefficients and a
  delay in dissipative boundary feedback.
\newblock {\em Journal of Mathematical Analysis and Applications},
  367(1):167--173, 2010.

\bibitem{rudy2019data}
Samuel Rudy, Alessandro Alla, Steven~L Brunton, and J~Nathan Kutz.
\newblock Data-driven identification of parametric partial differential
  equations.
\newblock {\em SIAM Journal on Applied Dynamical Systems}, 18(2):643--660,
  2019.

\bibitem{rudy2017data}
Samuel~H Rudy, Steven~L Brunton, Joshua~L Proctor, and J~Nathan Kutz.
\newblock Data-driven discovery of partial differential equations.
\newblock {\em Science Advances}, 3(4):e1602614, 2017.

\bibitem{sattar2020non}
Yahya Sattar and Samet Oymak.
\newblock Non-asymptotic and accurate learning of nonlinear dynamical systems.
\newblock {\em arXiv preprint arXiv:2002.08538}, 2020.

\bibitem{schaeffer2017learning}
Hayden Schaeffer.
\newblock Learning partial differential equations via data discovery and sparse
  optimization.
\newblock {\em Proceedings of the Royal Society A: Mathematical, Physical and
  Engineering Sciences}, 473(2197):20160446, 2017.

\bibitem{seymour1987exact}
Brian Seymour and Eric Varley.
\newblock Exact representations for acoustical waves when the sound speed
  varies in space and time.
\newblock {\em Studies in Applied Mathematics}, 76(1):1--35, 1987.

\bibitem{simchowitz2018learning}
Max Simchowitz, Horia Mania, Stephen Tu, Michael~I Jordan, and Benjamin Recht.
\newblock Learning without mixing: Towards a sharp analysis of linear system
  identification.
\newblock In {\em Conference On Learning Theory}, pages 439--473. PMLR, 2018.

\bibitem{vila2017bloch}
Javier Vila, Raj~Kumar Pal, Massimo Ruzzene, and Giuseppe Trainiti.
\newblock A bloch-based procedure for dispersion analysis of lattices with
  periodic time-varying properties.
\newblock {\em Journal of Sound and Vibration}, 406:363--377, 2017.

\bibitem{yang_stochastic_2020}
Zhenhuan Yang, Baojian Zhou, Yunwen Lei, and Yiming Ying.
\newblock Stochastic {Hard} {Thresholding} {Algorithms} for {AUC}
  {Maximization}.
\newblock In {\em 2020 {IEEE} {International} {Conference} on {Data} {Mining}
  ({ICDM})}, pages 741--750, Sorrento, Italy, November 2020. IEEE.

\bibitem{yilun_chen_sparse_2009}
{Yilun Chen}, Yuantao Gu, and Alfred~O. Hero.
\newblock Sparse {LMS} for system identification.
\newblock In {\em 2009 {IEEE} {International} {Conference} on {Acoustics},
  {Speech} and {Signal} {Processing}}, pages 3125--3128, Taipei, Taiwan, April
  2009. IEEE.

\bibitem{yuantao_gu_l_0_2009}
{Yuantao Gu}, {Jian Jin}, and {Shunliang Mei}.
\newblock $\ell_0$ norm constraint {LMS} algorithm for sparse system
  identification.
\newblock {\em IEEE Signal Processing Letters}, 16(9):774--777, September 2009.

\bibitem{zhai_tracking_2019}
Tingting Zhai, Frederic Koriche, Hao Wang, and Yang Gao.
\newblock Tracking {Sparse} {Linear} {Classifiers}.
\newblock {\em IEEE Transactions on Neural Networks and Learning Systems},
  30(7):2079--2092, July 2019.

\bibitem{zhang2019convergence}
Linan Zhang and Hayden Schaeffer.
\newblock On the convergence of the {SIND}y algorithm.
\newblock {\em Multiscale Modeling \& Simulation}, 17(3):948--972, 2019.

\bibitem{zinkevich2003online}
Martin Zinkevich.
\newblock Online convex programming and generalized infinitesimal gradient
  ascent.
\newblock In {\em Proceedings of the 20th international conference on machine
  learning (icml-03)}, pages 928--936, 2003.

\end{thebibliography}
